\documentclass[11pt]{article}
\usepackage{multicol}
\usepackage{amsfonts}
\usepackage{latexsym,amssymb,amsmath,graphics,cite}
\usepackage{fancyhdr}
\usepackage{lipsum}
\usepackage{lineno}
\usepackage{tikz}
\usepackage{CJK}
\usepackage{graphicx}
\usepackage{multirow}
\usepackage{float}
\usepackage{caption}
\usepackage{diagbox}
\usepackage{multirow}
\usepackage{setspace}
\usepackage{appendix}
\usepackage{tabu}
\begin{document}
\newcommand{\qed}{\hphantom{.}\hfill $\Box$\medbreak}
\newcommand{\Proof}{\noindent{\bf Proof \ }}
\newcommand{\rnote}[1]{{\color{red}{\sf #1}}}

\newtheorem{theorem}{Theorem}[section]
\newtheorem{lemma}[theorem]{Lemma}
\newtheorem{proposition}[theorem]{Proposition}
\newtheorem{corollary}[theorem]{Corollary}
\newtheorem{remark}[theorem]{Remark}
\newtheorem{example}[theorem]{Example}
\newtheorem{definition}[theorem]{Definition}
\newtheorem{Construction}[theorem]{Construction}
\begin{CJK*}{GBK}{song}
\title{\large{\bf Generalized perfect difference families and their application to
variable-weight geometric orthogonal codes \footnote{Supported by IPFHB under Grant CXZZBS2021067 (X. Su),
the Key R$\&$D Program of Guangdong province under Grant 2020B030304002 and NSFC under Grant 11901210 (L. Wang), NSFC under Grant 11871019 (Z. Tian)}}}

\author{Xiaowei Su$^{a}$, Lidong Wang$^b$, Zihong Tian$^a$\thanks{corresponding author} \\
\small $^a$ School of Mathematical Sciences, Hebei Normal University, Shijiazhuang 050024, P. R. China\\
\small $^b$ School of Intelligence Policing, China People's Police University, Langfang 065000, P. R. China\\
{\footnotesize 1367973156@qq.com, lidongwang@aliyun.com, tianzh68@163.com}}

\date{}

\maketitle

\begin{abstract}
 Motivated by the application in geometric orthogonal codes (GOCs), Wang et al.  introduced  the concept of generalized perfect difference families (PDFs), and established the equivalence between GOCs and a certain type of generalized PDFs recently.  Based on the relationship, we  discuss the existence problem of generalized $(n\times m,K,1)$-PDFs in this paper. By using some auxiliary designs such as semi-perfect group divisible designs and  several recursive constructions, we prove that a generalized $(n\times m, \{3,4\}, 1)$-PDF exists if and only if $nm\equiv1\pmod{6}$. The existence of a generalized $(n\times m, \{3,4,5\}, 1)$-PDF is  also completely solved possibly except for a few values.  As a consequence, some variable-weight perfect $(n\times m,K,1)$-GOCs are obtained.\vspace{0.2cm}

\end{abstract}

\section{Introduction}
Let $N $ and $M$ be two sets of integers containing $0$ such that if $a\in N$ (or $a\in M$) then $-a\in N$ (or $-a\in M$). Let  $K$ be a set of positive integers. For any $k$-subset $B$ of $N \times M$, where $k\in K$, define the multi-set
$$\Delta B :=\{(x_{1}-x_{2},y_{1}-y_{2}): (x_{1},y_{1}), (x_{2},y_{2})\in B, (x_{1},y_{1})\neq (x_{2},y_{2})\}$$
as the \emph{list of differences} from $B$. A \emph{generalized $(N \times M,K,1)$-perfect difference packing} (briefly PDP) $\mathcal{B}$ is a collection of subsets (called \emph{base~blocks}) of $N \times M$ of sizes from $K$ such that the multi-set $\Delta\mathcal{B}:=\bigcup_{B\in\mathcal{B}}\Delta B$ covers every element of $N \times M$ at most once. $(N \times M)\setminus \Delta\mathcal{B}$ is called the \emph{difference leave} or \emph{leave} of this packing. If the leave consists only of $(0,0)$, $\mathcal{B}$ is called a \emph{generalized $(N \times M,K,1)$-perfect difference family} (briefly PDF). As usual, a generalized $(N \times M,\{k\},1)$-PDP (or -PDF) will be abbreviated to a generalized $(N \times M,k,1)$-PDP (or -PDF).

\begin{example}\label{Example2}
The following $8$ base blocks form a generalized $(N \times M,3,1)$-PDF, where $N =\{0,\pm a,\pm b, \pm c\}$, $M=\{0,\pm x,\pm y, \pm z\}$, $a, b,c$ and $x,y,z$ are integers   satisfying that $a+b=c$, $x+y=z$.
\begin{center}
\begin{tabular}{llllll}
$\{(0, 0), (a, 0), (c, z)\}$,
$\{(0, x), (a, z), (c, 0)\}$,
$\{(0, z), (a, 0), (c, x)\}$,\\
$\{(0, z), (a, y), (c, 0)\}$,
$\{(0, 0), (a, z), (c, y)\}$,
$\{(0, 0), (0, z), (a, x)\}$,\\
$\{(0, 0), (0, y), (b, y)\}$,
$\{(0, 0), (0, x), (c, x)\}$.
\end{tabular}
\end{center}
\end{example}

\begin{example}\label{5x5}
The following $3$ base blocks form a generalized $(N \times M,\{3,4\},1)$-PDF, where $N =\{0,\pm 1,\pm 2\}$, $M=\{0,\pm 1,\pm 2\}$.
\begin{center}
\begin{tabular}{llllll}
$\{(0,0),(0,2),(2,0)\},$
$\{(0,1),(1,2),(2,0)\},$
$\{(0,0),(1,2),(2,1),(2,2)\}.$
\end{tabular}
\end{center}
\end{example}

For convenience, in what follows we always assume that $[a,b]=\{a,a+1,\dots,b\}$ for integers $a,b$ with $a<b$, and
$$[n]:=[-\frac{n-1}{2},\frac{n-1}{2}]$$
for any odd integer $n$, unless otherwise specified.
We usually write a generalized $(n\times m,K,1)$-PDP (or -PDF) instead of a generalized $([n]\times [m],K,1)$-PDP (or -PDF). Example \ref{5x5} provides a generalized $(5\times  5,\{3,4\},1)$-PDF.

When $m=1$, a generalized $(n\times m,K,1)$-PDP (or -PDF) is simply written as an $(n,K,1)$-PDP (or -PDF), and the second coordinates of all elements in the base blocks will be omitted.  The existence problem of $(n,K,1)$-PDPs was first introduced by Bermond, Kotzig and Turgeon \cite{Bermond} in 1976 in connection with a problem of spacing movable antennas in radioastronomy, who used a different terminology named {\em perfect systems of difference sets} (PSDSs). An $(e,\{k_1,k_2,\ldots,k_e\},c)$-PSDS is a family $S=\{S_1,S_2,\ldots,S_e\}$, where $S_i=\{s_{i1},s_{i2},\ldots,s_{ik_i}\}$, $0\leq s_{i1}<s_{i2}<\ldots<s_{ik_i}$ and all $s_{ij}$ are integers, such that
$\Delta S:=\bigcup_{1\leq i \leq e}\{s_{il}-s_{ij}:1\leq j<l\leq k_i\}=\{c,c+1,\ldots,c-1+\sum_{1\leq i \leq e}{k_i\choose2}\}$. Such a system is trivial if $k_i=2$ for at least one $i$.
An $(e,\{k_1,k_2,\ldots,k_e\},c)$-PSDS is \emph{regular} if $k_1=k_2=\cdots=k_e=k$. As usual,
a regular $(e, \{k\}, c)$-PSDS is abbreviated to $(e,k,c)$-PSDS, which can be regarded as an $(ek(k-1)+1,k,1)$-PDP with leave $[2c-1]$. Obviously, an $(n,k,1)$-PDF
is a regular $((n-1)/k(k-1),k,1)$-PSDS.

Bermond et al. \cite{Bermond} proved that perfect difference families do not exist for any $k\geq6$. For $k=3$, the existence problem has been completely settled.
However, the existence problem for $k=4$ is far from being settled. The existence results for
$k=5$ are even more scarce. Here, we list some known results for $(n,k,1)$-PDFs.

\begin{lemma}{\rm (Abel et al. \cite{Difference-families}, Beth et al. \cite{Beth}, Ge et al. \cite{GeG})}\label{n,3,1}\\
$(1)$~There exists an $(n,3,1)$-PDF if and only if $n\equiv1,7\pmod{24};$\\
$(2)$~There exists a $(12t+1,4,1)$-PDF for every $t\leq1000$ with two definite exceptions of $t=2,3;$\\
$(3)$~There exists a $(20t+1,5,1)$-PDF for $t=6,8,10;$\\
$(4)$~There are no $(n,5,1)$-PDFs for $n\equiv21\pmod{40}$~or~$n=41,81.$
\end{lemma}

In 1983, Laufer et al. \cite{Laufer} proved a conjecture of Erd\"os, that is for every positive integer $t$, except for a finite number, there is a non-trivial PSDS whose differences are the first $t$ positive integers, i.e., there is a $(2t+1,K,1)$-PDF where $k\geq3$ for all $k\in K$. $(n,K,1)$-PDFs are interesting not only in their own right, but also in their relation to many other combinatorial configurations such as semi-cyclic holey group divisible designs \cite{Feng} and optimal variable-weight optical orthogonal codes \cite{Jiang, WuD, yang}. We quote some results for later use. Note that an $(n,\{k,s^{*}\},1)$-PDF is an $(n,\{k,s\},1)$-PDF with exactly one base block of size $s$, and the number of base blocks of size $k$ is greater than zero. An $(n,\{k,s_1^{*},s_2^{*}\},1)$-PDF is an $(n,\{k,s_1,s_2\},1)$-PDF with exactly one base block of size $s_i$, $i=1,2$ and the number of base blocks of size $k$ is greater than zero.

\begin{lemma}{\rm (Wu et al. \cite{WuD}, Chen \cite{ChenZ}, Feng et al. \cite{Feng})}\label{n-345}\\
$(1)$~There exists an $(n,\{3,4^{*}\},1)$-PDF if and only if $n\equiv1\pmod{6}$ and $n\geq19;$\\
$(2)$~There exists an $(n,\{3,5^{*}\},1)$-PDF if and only if $n\equiv9,15\pmod{24}$ and $n\geq33;$\\
$(3)$~There exists an $(n,\{3,6^{*}\},1)$-PDF if and only if $n\equiv1\pmod{6}$ and $n\geq43;$\\
$(4)$~There exists an $(n,\{3,7^{*}\},1)$-PDF if and only if $n\equiv1,7\pmod{24}$ and $n\geq73;$\\
$(5)$~There exists an $(n,\{3,4^{*},5^{*}\},1)$-PDF if and only if $n\equiv3\pmod{6}$ and $n\geq39;$\\
$(6)$~There exists an $(n,\{3,4,5\},1)$-PDF for $n\equiv5\pmod{6}$ and $n\geq59$.
\end{lemma}

When $m>1$, generalized $(n\times m,k,1)$-PDFs were first introduced by Wang et al. \cite{WLD}, which was motivated by the application in constructions of geometric orthogonal codes. They gave a complete solution to the existence of a generalized $(n\times m,3,1)$-PDF.

\begin{lemma}{\rm (Wang et al. \cite{WLD})}\label{nxm,3,1}
There exists a generalized $(n\times m,3,1)$-PDF if and only if $nm\equiv1\pmod{6}$ and $n, m\equiv1,7,17,23\pmod{24}$.
\end{lemma}

This paper continues to discuss the existence of generalized $(n\times m,K,1)$-PDFs with $|K|\geq2$. The reason is that, on the one hand, the existence problem of generalized $(n\times m,K,1)$-PDFs can be regarded as a generalization of the conjecture of Erd\"os (considering the elements in the base block from an integer to a pair of integers), on the other hand, generalized $(n\times m,K,1)$-PDFs can be used to construct variable-weight geometric orthogonal codes. As the main results of this paper, we are to prove the following theorems.

\begin{theorem}\label{result34}
There is a generalized $(n\times m, \{3,4\}, 1)$-PDF if and only if $nm\equiv1\pmod{6}$.
\end{theorem}

\begin{theorem}\label{result345}
There is a generalized $(n\times m, \{3,4,5\}, 1)$-PDF if and only if $n, m\equiv1\pmod{2}, n,m\neq3$, except for $\{n,m\}\in\{\{1,d\},\{5,7\},\{5,9\}\}$, where~$d\in \{5,9,11,15,17,21,23,$ $27,29,35,41,47,53\}$ and possibly except for $\{n,m\}\in\{\{5,13\},\{5,45\},\{7,23\},\{7,29\},\{7,35\},$ $\{9,35\},\{11,19\},\{11,27\},\{13,17\},\{13,21\},\{13,23\},\{13,27\},\{13,29\},\{13,35\},\{15,17\},\{15,\\ 21\},\{15, 27\},\{17,21\},\{17, 27\},\{21, 21\},\{21, 27\},\{23, 45\},\{27, 27\},\{29, 45\},\{35,45\}\}$.
\end{theorem}

The rest of the paper is organized as follows. In Section 2, we introduce  some auxiliary designs. In Section 3, we present several recursive constructions for generalized perfect difference packings. In Section 4, we show that the necessary conditions for the existence of a generalized $(n\times m,\{3,4\},1)$-PDF  are  also sufficient. The existence of a generalized $(n\times m,\{3,4,5\},1)$-PDF is almost completely given except for a few possible values in Section 5. In Section 6, we consider the application of generalized perfect difference families to variable-weight geometric orthogonal codes. Conclusions and discussions are given in Section 7.

\section{Auxiliary designs}
In this section, we introduce some auxiliary designs, which will play important roles in the
recursive constructions for generalized perfect difference families.

\subsection{Semi-perfect group divisible designs}
Let $K$ be a set of positive integers. A \emph{group divisible design}, denoted by $K$-GDD, is a triple
$(X, \mathcal{G}, \mathcal{B})$ satisfying that: (1) $X$ is a finite set of  points;
(2) $\mathcal{G}$ is a partition of $X$ into  subsets (called \emph{groups});
(3) $\mathcal{B}$ is a collection of subsets (called \emph{blocks}) of $X$, each
block has cardinality from  $K$, such that every 2-subset of $X$ is either contained in exactly one block or in exactly one group, but not in both. The
multi-set $T = \{|G| : G \in \mathcal{G}\}$ is called the \emph{type} of the $K$-GDD. If $\mathcal{G}$ contains $u_i$ groups of size $g_i$ for $1\leq i\leq r$, then we also denote the type by $g_1^{u_1} g_2^{u_2}\dots g_r^{u_r}$.  If $K=\{k\}$, we write a $\{k\}$-GDD as a $k$-GDD.

Let $S$ be a set of $n$ points and $X=S\times[m]$, $\mathcal{G}=\{\{i\}\times [m]$: $i\in S\}$.
Let $\mathcal{F}$ be a collection of subsets (called \emph{base blocks}) of $S\times[m]$. For $i,j\in S$ and $B\in{\cal F}$, define a multi-set $\Delta_{ij}(B)=\{x-y: (i,x),(j,y)\in B, (i,x)\neq(j,y)\}$, and a multi-set $\Delta_{ij}({\cal F})=\bigcup_{B\in{\cal F}}\Delta_{ij}(B)$. If for any $i,j\in S$,
$$\Delta_{ij}({\cal F})=\left\{
\begin{array}{lll}
[m], & \ \ i\neq j,\\
\emptyset, & \ \ i=j,\\
\end{array}
\right.
$$
then a $K$-GDD of type $m^n$ with the point set $X'=S\times [0,m-1]$ together with the group set ${\cal {G}}' =\{\{i\} \times [0,m-1]: i\in S\}$ can be generated from ${\cal F}$, where $K=\{|B|:B\in{\cal F}\}$. The required blocks are
obtained by developing all base blocks of ${\cal F}$ by successively
adding $1$ to the second component of each point of these base
blocks modulo $m$. We call $\mathcal{F}$ a \emph{semi-perfect group divisible design}, denoted by $K$-SPGDD of type $m^n$.
As usual, a $\{k\}$-SPGDD will be abbreviated to a $k$-SPGDD.

\begin{example}\label{9SPGDD}
All base blocks of a $\{3,4\}$-SPGDD of type $9^4$ on $[0,3]\times [9]$ are listed below.\vspace{0.2cm}

\hspace{0.5cm}$\begin{array}{lll}
\{(0, 4), (1, 0), (2, 1)\},~~~ & \{(0, 0), (1, 4), (2, 0)\},~~~ & \{(1, 1), (2, 0), (3, 4)\},\\
\{(0, 4), (1, 1), (3, 2)\}, & \{(0, 0), (2, 1), (3, 4)\}, & \{(0, 0), (1, 0), (3, 0)\},\\
\{(0, 1), (2, 4), (3, 0)\}, & \{(1, 0), (2, 4), (3, 2)\}, & \{(0, 3), (1, 4), (2, 1), (3, 0)\},\\
\{(0, 1), (1, 4), (3, 3)\}, & \{(0, 0), (1, 2), (2, 4)\}, & \{(0, 1), (1, 0), (2, 3), (3, 4)\},\\
\{(1, 3), (2, 3), (3, 0)\}, & \{(0, 1), (2, 0), (3, 2)\}, & \{(0, 4), (1, 2), (2, 0), (3, 0)\}.
\end{array}$
\end{example}

\begin{example}\label{11SPGDD}
All base blocks of a $\{3,4\}$-SPGDD of type $11^4$ on $[0,3]\times [11]$ are listed below.\vspace{0.2cm}

\hspace{0.5cm}$\begin{array}{lll}
\{(0, 0), (1, 1), (3, 1)\},~~~ & \{(0, 0), (1, 3), (3, 5)\},~~~ & \{(0, 0), (1, 0), (2, 0)\},\\
\{(0, 5), (1, 3), (3, 0)\},~ & \{(0, 0), (2, 5), (3, 4)\},~ & \{(1, 0), (2, 3), (3, 1)\},\\
\{(0, 1), (1, 5), (2, 0)\},~ & \{(0, 5), (2, 0), (3, 2)\},~ & \{(1, 1), (2, 0), (3, 4)\},\\
\{(1, 5), (2, 3), (3, 0)\},~ & \{(1, 0), (2, 2), (3, 5)\},~ & \{(0, 2), (1, 4), (2, 0), (3, 0)\},\\
\{(0, 3), (2, 0), (3, 5)\},~ & \{(0, 5), (1, 0), (2, 1)\},~ & \{(0, 4), (1, 1), (2, 5), (3, 0)\},\\
\{(0, 4), (1, 0), (3, 4)\},~ & \{(0, 1), (2, 4), (3, 0)\},~ & \{(0, 0), (1, 5), (2, 2), (3, 3)\},\\
\{(0, 1), (1, 0), (2, 5)\}.~& &
\end{array}$
\end{example}

SPGDDs are closely related to \emph{perfect difference matrices} (PDMs).  A $k\times m$ matrix $D=(d_{ij})$ with entries from $[m]$ is called a \emph{perfect difference matrix}, denoted by  PDM$(k,m)$, if for all $1\leq s<t\leq k$, the list of differences $\{d_{sj}-d_{tj}:1\leq j\leq m\}$ contains each element of $[m]$ exactly once.

We consider a  $k$-SPGDD  of type $m^k$ on $S\times [m]$, $S=\{s_0,s_1,\dots,s_{k-1}\}$. It is easy to know that such an SPGDD contains exactly $m$ base blocks. For each base block $B=\{(s_0,x_0),(s_1,x_1),\dots,(s_{k-1},x_{k-1})\}$, if we regard $(x_0,x_1,\dots,x_{k-1})^T$ as a column of a $k\times m$ matrix $D$, then it is readily checked that $D$ is a PDM($k,m$). This means that the concept of a semi-perfect group divisible design is in fact a generalization of  a perfect difference matrix.

\begin{proposition}
A $k$-SPGDD  of type $m^k$ is equivalent to a PDM$(k,m)$.
\end{proposition}
\newpage
Perfect difference matrices were investigated by Ge et al. \cite{GeG} for the constructions of perfect difference families and additive sequences of permutations. Note that the PDM$(k,m)$ defined in \cite{GeG} also require that the entries of each row  comprise all the elements of $[m]$, then a PDM$(k,m)$ defined in this paper is equivalent to a PDM$(k-1,m)$ defined in \cite{GeG}.

\begin{lemma}{\rm(Ge et al. \cite{GeG}, Wang et al. \cite{WLD})}\label{345PDM}\\
$(1)$~There exists a PDM$(3,m)$ for any odd integer $m\geq3;$\\
$(2)$~There exists a PDM$(4,m)$~for any odd integer $3<m<200$ expect for $m=9,11$ and possibly for $m=59;$\\
$(3)$~There exists a PDM$(5,m)$ for any $m=5^{a_1}121^{a_2}161^{a_3}201^{a_4}$, where $a_1,a_2,a_3,a_4$ are non-negative integers, not all equal to 0.
\end{lemma}

\subsection{Semi-perfect modified group divisible designs}

A {\em modified group divisible design} (MGDD) is a quadruple $(X,{\cal H},{\cal G},{\cal B})$ which satisfies the following properties: $(1)$ $X$ is a finite set of $mn$ {\em points}; $(2)$ ${\cal G}$ is a partition of $X$ into $n$ subsets (called {\em groups}), each of size $m$; $(3)$ ${\cal H}$ is another partition of $X$ into $m$ subsets (called {\em holes}), each of size $n$ such that $|H\cap G|=1$ for any $H\in{\cal H}$ and $G\in{\cal G}$; $(4)$ $\cal B$ is a set of subsets (called {\em blocks}) of $X$, each
block has cardinality from a set of
positive integers $K$, such that no block contains two distinct points of any group or any hole, but any other pair of distinct points of $X$ occurs in exactly one block of $\cal B$.
Such a design is denoted by a $K$-MGDD of type $m^n$. If $K=\{k\}$, we write a $\{k\}$-MGDD as a $k$-MGDD.

Assaf \cite{3-MGDD} first introduced the concept of MGDDs and gave the existence of 3-MGDDs. The existence problem of 4-MGDDs was  completely resolved in \cite{caoh}. There are also some results  on 5-MGDDs (See \cite{MGDD,5-MGDD}). Here,  we only quote the following results for later use.

\begin{lemma}{\rm (Abel et al. \cite{MGDD}, Cao et al. \cite{caoh})}\label{MGDD}\\
$(1)$ There exists a $3$-MGDD of type $m^n$ if and only if $n, m\geq3$, $(n-1)(m-1)\equiv0\pmod {2}$ and $n(n-1)m(m-1)\equiv0\pmod {6};$\\
$(2)$ There exists a $4$-MGDD of type $m^n$ if and only if $n, m\geq4$ and $(n-1)(m-1)\equiv0\pmod {3}$ except for $(n,m)=(4,6);$\\
$(3)$ There exists a $5$-MGDD of type $5^5$.
\end{lemma}

Let $S$ be a set of $n$ points and $X=S\times[m]$, $\mathcal{G}=\{\{i\}\times [m]$: $i\in S\}$ and $\mathcal{H}=\{S\times\{x\}:x\in[m]\}$.
Let $\mathcal{F}$ be a collection of subsets (called \emph{base blocks}) of $S\times[m]$. For $i,j\in S$ and $B\in{\cal F}$, define a multi-set $\Delta_{ij}(B)=\{x-y: (i,x),(j,y)\in B, (i,x)\neq(j,y)\}$, and a multi-set $\Delta_{ij}({\cal F})=\bigcup_{B\in{\cal F}}\Delta_{ij}(B)$. If for any $i,j\in S$,
$$\Delta_{ij}({\cal F})=\left\{
\begin{array}{lll}
[m]\setminus\{0\}, & \ \ i\neq j,\\
\emptyset, & \ \ i=j,\\
\end{array}
\right.
$$
then a $K$-MGDD of type $m^n$ with the point set $X'=S\times [0,m-1]$ together with the group set ${\cal {G}}' =\{\{i\} \times [0,m-1]: i\in S\}$ and the hole set ${\cal H}'=\{S\times \{x\}: x\in[0,m-1]\}$ can be generated from ${\cal F}$, where $K=\{|B|:B\in{\cal F}\}$. The required blocks are
obtained by developing all base blocks of ${\cal F}$ by successively
adding $1$ to the second component of each point of these base
blocks modulo $m$. We call $\mathcal{F}$ a \emph{semi-perfect modified group divisible design}, denoted by $K$-SPMGDD of type $m^n$.
Clearly, if $\mathcal{F}$ is a $K$-SPMGDD of type $m^n$ on $S\times [m]$, then $\mathcal{F}\cup\{S\times\{0\}\}$ forms a $(K\cup\{n\})$-SPGDD of type $m^n$.

\section{Construction methods}
In this section, we shall establish some recursive constructions for generalized perfect difference packings.

Let $H_1$ and $H_2$ be two sets of integers. For any positive integers $r_1$ and $r_2$, define
$$H_1^{r_1}\times H_2^{r_2}:=\{(r_1 x_1, r_2 x_2):(x_1,x_2)\in H_1\times H_2\}.$$
When $H_i=[h_i]$ for some odd positive integer $h_i$, we write $H_i^{r_i}$ as $[h_i]^{r_i}$. When $r_i=1$, we omit $r_i$ from the notation. A generalized $(N \times M,K,1)$-PDF can be viewed as a generalized $(N \times M,K,1)$-PDP with leave $[1]^{r_1}\times [1]^{r_2}$ for any positive integers $r_1$ and $r_2$. Further, if there exists a generalized $(N \times M,K,1)$-PDP with leave $H_1^{r_1}\times H_2^{r_2}$, then $H_i$ contains $0$ and if $a\in H_i$, then $-a\in H_i$, where $i=1,2$.
 The following construction is a specialization of Construction 4.1 in \cite{WLD}.

\begin{Construction}\label{budong}
Suppose that there exist
\begin{itemize}
\item[$(1)$] a generalized $(n\times m,K,1)$-PDP with leave $H_1^{r_1}\times H_2^{r_2}$, and
\item[$(2)$] a generalized $(H_1\times H_2,L,1)$-PDP with leave $T_1^{s_1}\times T_2^{s_2}$,
\end{itemize}
then there exists a generalized $(n\times m,K\cup L,1)$-PDP with leave $T_1^{r_1s_1}\times T_2^{r_2s_2}$.
Moreover, if the generalized $(H_1\times H_2,L,1)$-PDP is a PDF, then the resulting PDP is also a PDF.
\end{Construction}

{\Proof}Let $\mathcal{F}$ be the given generalized $(n\times m,K,1)$-PDP with leave $H_1^{r_1}\times H_2^{r_2}$ on $[n]\times[m]$.
Let $\mathcal{A}$ be a generalized $(H_1\times H_2,L,1)$-PDP with leave $T_1^{s_1}\times T_2^{s_2}$.
For each base block $A=\{(a_1,b_1),(a_2,b_2),\dots,(a_l,$ $b_l)\}  \in\mathcal{A}$ with $l\in L$, construct a set $B_A=\{(r_1a_1,r_2b_1),(r_1a_2,r_2b_2),$ $\dots,(r_1a_l,r_2b_l)\}$. Let $\mathcal{B}=\{B_A: A\in \mathcal{A}\}$. Then $\mathcal{F}\cup\mathcal{B}$ forms a generalized $(n\times m,K\cup L,1)$-PDP with leave $T_1^{r_1s_1}\times T_2^{r_2s_2}$. \qed

\begin{Construction}\label{N-NxM-PDP}
Suppose that there exists an $(n,K,1)$-PDP with leave $H^r$. If
there exists an $L$-SPGDD of type $m^k$ for each $k\in K$, then
there exists a generalized $(n\times m, L,1)$-PDP with leave $H^r\times [m]$.
\end{Construction}
{\Proof}Let $\mathcal{A}$ be the given $(n,K,1)$-PDP  with leave $H^r$ on $[n]$.
 For any $A\in\mathcal{A}$, we can construct by assumption an $L$-SPGDD of type $m^k$ on $A\times [m]$ with group set $\{\{i\}\times [m]: i\in A\}$ and base block set $\mathcal{B}_A$, where $k\in K$.
 It is not difficult to find $\Delta  {\mathcal{B}}_A=\Delta A\times [m]$.
   Define
 $$\mathcal{B}=\bigcup \limits_{A\in \mathcal{A}}\mathcal{B}_{A},$$
 then
 $\Delta \mathcal{B}=\bigcup _{A\in \mathcal{A}} \Delta {\mathcal{B}}_A =([n]\setminus H^r)\times [m].$
 Therefore, $\mathcal{B}$ forms a generalized  $(n\times m, L,1)$-PDP with leave $H^r\times [m]$.\qed

\noindent\textbf{Remark}~~If there is an $(n,K,1)$-PDF in Construction \ref{N-NxM-PDP}, then the leave of the resulting PDP is $\{0\}\times [m]$.

\begin{Construction}\label{MGDD-SPMGDD}
Suppose there exists an $(m,K,1)$-PDF. If there exists an $L$-MGDD of type $k^h$ for each $k\in K$, then there exists an $L$-SPMGDD of type $m^h$.
\end{Construction}
{\Proof}Let $\mathcal{A}$ be the given   $(m,K,1)$-PDF on $[m]$.
 For each $A=\{x_1,x_2,\dots,x_k\}\in\mathcal{A}$, let $\mathcal{B}_A$ be an  $L$-MGDD of type $k^h$ on $[0,h-1]\times A$ with  group set $\{\{i\}\times A: i\in [0,h-1]\}$, where $k\in K$. For each $B\in\mathcal{B}_A$, let $\Delta_{ij} B=\{x-y: (i,x),(j,y)\in B\}$, where $i,j\in [0,h-1]$ and $i\neq j$.
Denote
$$\mathcal{B}=\bigcup \limits_{A\in \mathcal{A}}\mathcal{B}_A.$$
From the properties of MGDD, we get
$$\Delta_{ij} \mathcal{B} =\bigcup \limits_{A\in \mathcal{A}}\Delta_{ij} \mathcal{B}_A=\bigcup \limits_{A\in \mathcal{A}}\bigcup \limits_{B\in \mathcal{B}_A}\Delta_{ij} B=[m]\setminus\{0\},~i\neq j\in [0,h-1].$$
It is straightforward to check that $\mathcal{B}$ constitutes an $L$-SPMGDD of type $m^h$ on $[0,h-1]\times[m]$.\qed

\begin{Construction}\label{nxm-nxmv}
Suppose that there exists a generalized $(n\times m,K,1)$-PDF. If there exists an $L$-SPGDD of type $v^k$ for each $k\in K$, and a  $(v,L',1)$-PDF, then there exists a generalized $(n\times(mv),L\cup L',1)$-PDF.
\end{Construction}
{\Proof} Let $\mathcal{A}$ be the given  generalized $(n\times m,K,1)$-PDF on $[n]\times[m]$.
For each $A=\{(x_0,y_0),(x_1,y_1),\dots,$ $(x_{k-1},y_{k-1})\}\in\mathcal{A}$ with $k\in K$, we can construct by assumption an $L$-SPGDD of type $v^k$ on $A\times [v]$  with group set $\{\{(x_i,y_i)\}\times [v]: 0\leq i \leq k-1\}$ and base block set $\mathcal{B}_A.$  For each $B=\{(x_{i_1},y_{i_1},z_{i_1}),(x_{i_2},y_{i_2},z_{i_2}),\dots,(x_{i_l},y_{i_l},z_{i_l})\}\in\mathcal{B}_A$ with $l\in L$, define
$$B_A=\{(x_{i_1},y_{i_1}+mz_{i_1}),(x_{i_2},y_{i_2}+mz_{i_2}),\dots,(x_{i_l},y_{i_l}+mz_{i_l})\}.$$
Let
$$\mathcal{C}=\bigcup \limits_{A\in \mathcal{A}}(\bigcup \limits_{B\in \mathcal{B}}B_A).$$
Then
$$\Delta\mathcal{C}=\bigcup \limits_{A\in \mathcal{A}}(\bigcup \limits_{B\in \mathcal{B}}\Delta B_A)
=\bigcup \limits_{A\in \mathcal{A}}\{(\alpha,\beta+m[v]):(\alpha,\beta)\in\Delta A\}=[n]\times [mv]\setminus\big\{\{0\}\times m[v]\big\}.$$

Let $\mathcal{D}$ be a  $(v,L',1)$-PDF. For each base block $D=\{a_1,a_2,\dots,a_s\}\in\mathcal{D}$ with $s\in L'$, construct a set $F_D=\{(0,ma_1),(0,ma_2),\dots,(0,ma_s)\}$. Let $\mathcal{F}=\{F_D:D\in \mathcal{D}\}$, then $\Delta \mathcal{F}=\{\{0\}\times m[v]\}\setminus\{(0,0)\}$.
\vspace{0.2cm}

It is not difficult to check that $\mathcal{C}\cup\mathcal{F}$ forms a generalized $(n\times(mv),L\cup L',1)$-PDF on  $[n]\times[mv]$.
\qed

\begin{corollary}\label{nxmv}
Suppose that there exists a generalized $(n\times m,\{3,4,5\},1)$-PDF and a $(v,$ $\{3,4,5\},1)$-PDF, then there is a generalized $(n\times(mv),\{3,4,5\},1)$-PDF.
\end{corollary}
{\Proof}There exists a $\{3,4,5\}$-MGDD of type $g^k$ for any $g,k\in \{3,4,5\}$ by Lemma \ref{MGDD}. Applying Construction \ref{MGDD-SPMGDD} with a  $(v,\{3,4,5\},1)$-PDF, we get a $\{3,4,5\}$-SPMGDD of type $v^k$ on $A\times [v]$ for each $k\in \{3,4,5\}$, where $A$ is a set of $k$ points. Then we get a  $\{3,4,5\}$-SPGDD of type $v^k$ by adding the base block $\{A\times\{0\}\}$. Take $K=L=L'=\{3,4,5\}$ in Construction \ref{nxm-nxmv}, we get the conclusion immediately.\qed

{\em Langford sequence} is an important combinatorial structure,  which has been used to construct  perfect difference packings by Chee et al. \cite{cklh}.

A {\em Langford sequence of order $n$ and defect $d$}, $n>d$, is a partition of $[1,2n]$ into a collection of ordered pairs $(a_i,b_i)$ such that $\{b_i-a_i : 1 \leq i \leq n\}=[d,d+n- 1]$.

\begin{lemma}{\rm (Simpson \cite{LD})}\label{LDjieguo}
A Langford sequence of order $n$ and defect $d$ exists if and only if
~$(1)$~$n\geq2d-1$, and
~$(2)$~$n\equiv0,1\pmod{4}$ for $d$  odd, or $n\equiv0,3\pmod{4}$ for $d$  even.
\end{lemma}

\begin{Construction}{\rm (Chee et al. \cite{cklh})}\label{LD-PDP}
If there exists a Langford sequence of order $n$ and defect $d$, then there exists a $(6n+2d-1,3,1)$-PDP with leave $[2d-1]$.
\end{Construction}

\section{Generalized $(n\times m, \{3,4\},1)$-PDFs}
In this section, we will discuss the existence problem of generalized $(n\times m,\{3,4\},1)$-PDFs. At first, we consider the case $m=1$.

\begin{lemma}\label{yiwei34}
There exists an $(n,\{3,4\},1)$-PDF if and only if $n\equiv1\pmod{6}$.
\end{lemma}
{\Proof}Let $\mathcal{A}$ be an $(n,\{3,4\},1)$-PDF on $[n]$, then
$$\Delta \mathcal{A}=\bigcup \limits_{B\in \mathcal{A}}\Delta B=[n]\setminus\{0\}.$$
Thus, there are $n-1$ different elements in $\Delta \mathcal{A}$. Denote the number of base blocks of sizes $3, 4$ in $\mathcal{A}$ by $x, y$ respectively. By the definition of a perfect difference family, we have $6x+12y=n-1$, which implies $n\equiv1\pmod{6}$. This gives the necessary condition.

It is easy to check that $\{\{0,1,3\}\}$ constitutes a   $(7,3,1)$-PDF, and $\{\{0,2,5,6\}\}$ constitutes a   $(13,4,1)$-PDF. Therefore, the sufficient condition is given by Lemma \ref{n-345}.\qed

\begin{lemma}\label{SPGDD34}
There exists a $\{3,4\}$-SPGDD of type $m^4$ for $m\equiv1\pmod{6}$.
\end{lemma}
{\Proof} Applying Construction \ref{MGDD-SPMGDD} with an $(m,\{3,4\},1)$-PDF in Lemma \ref{yiwei34} and a \{3,4\}-MGDD of type $k^4$ for $k\in\{3,4\}$ in Lemma \ref{MGDD}, we get  a \{3,4\}-SPMGDD of type $m^4$ on $A\times [m]$ for $m\equiv1\pmod{6}$  with base block set $\mathcal{F}$, where $A$ is a set of $4$ points. Furthermore, let
$$\mathcal{C}=\mathcal{F}\cup\{A\times\{0\}\},$$
then $\mathcal{C}$ forms a $\{3,4\}$-SPGDD of type $m^4$.\qed

\begin{lemma}\label{34CY}
There exists a generalized $(n\times m,\{3,4\},1)$-PDF for $n,m\equiv1\pmod{6}$.
\end{lemma}
{\Proof}When $n,m\equiv1\pmod{6}$, there is an $(n,\{3,4\},1)$-PDF and an   $(m,\{3,4\},1)$-PDF from Lemma \ref{yiwei34}. By Lemma \ref{SPGDD34}, there is a $\{3,4\}$-SPGDD of type $m^4$. By Lemma \ref{345PDM}, there is a PDM(3,$m$), i.e., a 3-SPGDD of type $m^3$. Therefore, we get a generalized $(n\times m, \{3,4\},1)$-PDF by  Construction \ref{budong} and Construction \ref{N-NxM-PDP}.\qed

A generalized $(n\times m,K,1)$-PDP with leave $[{h_1}]^{r_1}\times[{h_2}]^{r_2}$ plays an important role in constructing generalized $(n\times m,K,1)$-PDFs.
 Wang et al. gave the following results.
\newpage
\begin{lemma}{\rm (Wang et al. \cite{WLD})}\label{5,11,17,23}\\
$(1)$~When $n,m\equiv5,11\pmod{24}$, there exists a generalized $(n\times m,3,1)$-PDP with leave $[{h_1}]^{r_1}\times[{h_2}]^{r_2}$, where $\{h_1,h_2\}\in\{\{5,5\},\{5,11\},\{11,11\}\}$.\\
$(2)$~When $n\equiv5,11\pmod{24},m\equiv17,23\pmod{24}$, there exists a generalized $(n\times m,3,1)$-PDP with leave $[{h_1}]^{r_1}\times[{h_2}]^{r_2}$, where $\{h_1,h_2\}\in\{\{5,5\},\{5,11\},\{5,17\},\{11,11\},\{11,17\}\}$.
\end{lemma}

\begin{lemma}\label{3,4XCS}
There exists a generalized $(n\times m,\{3,4\},1)$-PDF for $\{n,m\}\in\{\{5,5\},\{5,11\},\{5,$ $17\},\{11,11\},\{11,17\}\}$.
\end{lemma}

\Proof For $(n,m)\in\{(5,11),(5,17)\}$, the base blocks of a generalized $(n\times m,\{3,4\},1)$-PDF are listed below.\vspace{0.2cm}

~$(n,m)=(5,11)$:
{\begin{quote}
$\{(0,0),(0,4),(2,5)\},~~\{(0,0),(1,1),(2,4)\},~~\{(0,2),(2,0),(2,5)\},$\\
$\{(0,0),(0,1),(1,5)\},~~\{(0,0),(0,2),(2,2)\},~~\{(0,3),(1,5),(2,0)\},$\\
$\{(0,5),(1,3),(2,0)\},~~\{(0,4),(1,4),(2,0),(2,3)\}.$
\end{quote}}

~$(n,m)=(5,17)$:
{\begin{quote}
$\{(0,0),(1,0),(2,8)\},~~\{(0,0),(0,8),(2,0)\},~~\{(0,0),(0,1),(2,6)\},~~\{(0,1),(1,8),(2,0)\},$\\
$\{(0,6),(1,1),(2,0)\},~~\{(0,7),(1,0),(2,2)\},~~\{(0,4),(0,7),(2,0)\},~~\{(0,0),(0,2),(1,6)\},$\\
$\{(0,0),(0,6),(2,7)\},~~\{(0,2),(2,0),(2,5)\},~~\{(0,0),(1,5),(2,2)\},~~\{(0,2),(0,6),(1,0)\},$\\
$\{(0,0),(0,7),(1,3),(2,4)\}.$
\end{quote}}
The generalized $(5\times5,\{3,4\},1)$-PDF see Example \ref{5x5}. For $(n,m)\in\{(11,11),(11,17)\}$, the base blocks of generalized $(n\times m,\{3,4\},1)$-PDFs see Appendix A. Note that the existence of a generalized $(n\times m,\{3,4\},1)$-PDF implies the existence of a generalized $(m\times n,\{3,4\},1)$-PDF. So, the conclusion holds.
\qed

\emph{Proof of Theorem} \ref{result34}

It is not difficult to verify that a generalized $(n\times m,K,1)$-PDP implies an $(nm,K,1)$-PDP with the same number of base blocks by a mapping
$\tau:(x,y)\mapsto x+yn,$  $(x,y)\in [n]\times[m],~x+yn\in [nm]$. Then the necessity is immediately obtained by Lemma \ref{yiwei34}. Next we give the proof of the sufficiency.

When $n,m\equiv1\pmod{6}$, the conclusion is given by Lemma \ref{34CY}. Now we only need to consider the case $n,m\equiv5\pmod{6}$. For $n,m\equiv17,23\pmod{24}$, there exists a generalized $(n\times m,3,1)$-PDF by Lemma \ref{nxm,3,1}. For $n,m\equiv5,11\pmod{24}$ or $n\equiv5,11\pmod{24}, m\equiv17,23\pmod{24}$, there is a generalized $(n\times m,3,1)$-PDP with leave $[{h_1}]^{r_1}\times[{h_2}]^{r_2}$ from Lemma \ref{5,11,17,23}, where $\{h_1,h_2\}\in\{\{5,5\},\{5,11\},\{5,17\},\{11,11\},\{11,17\}\}$. Further, we can get these generalized $(h_1\times h_2,\{3,4\},1)$-PDFs from Lemma \ref{3,4XCS}. Applying Construction \ref{budong}, we get a generalized $(n\times m,\{3,4\},1)$-PDF for $nm\equiv1\pmod{6}$.\qed\par

\section{Generalized $(n\times m, \{3,4,5\},1)$-PDFs}
In this section, we investigate the existence of generalized $(n\times m,\{3,4,5\},1)$-PDFs. We first consider the case $m=1$, then some generalized $(n\times m, \{3,4,5\},1)$-PDFs can be obtained from $(n, \{3,4,5\},1)$-PDFs by applying the constructions in Section $3$.  With the help of computer search, some generalized PDFs  and PDPs with small parameters are also obtained. In addition, a generalized $(n\times m,\{3,4\},1)$-PDF is obviously also a generalized $(n\times m,\{3,4,5\},1)$-PDF.

\subsection{$(n,\{3,4,5\},1)$-PDFs}
\begin{lemma}\label{necessary}
If there is an $(n,\{3,4,5\},1)$-PDF, then $n\equiv1\pmod{2}$ and $n\notin\{3,5,9,11,15,$ $17,23,29,35\}$.
\end{lemma}
\Proof Suppose $\mathcal{F}$ is an $(n,\{3,4,5\},1)$-PDF and the number of base blocks of sizes $3,4$ and $5$ in $\mathcal{F}$ is  $x,y$ and $z$, respectively. Considering the list of differences for any base block in $\mathcal{F}$, we have
\begin{equation}\label{nece}
6x+12y+20z=n-1.
\end{equation}

It is easy to see that Equation (\ref{nece}) has a non-negative integer solution $(x,y,z)$ only if  $n\equiv 2z+1\pmod{6}$. Therefore, there is no  $(n,\{3,4,5\},1)$-PDF for $n\equiv3\pmod{6}$ and $n<21$, or $n\equiv5\pmod{6}$ and $n<41$.
The conclusion  holds.
\qed

\begin{lemma}\label{yiwei}
There is an $(n,\{3,4,5\},1)$-PDF if and only if $n\equiv1\pmod{2}$ and  $n\notin \{3,5,9,11,$ $15,17,21,23,27,29,35,41,47,53\}$.
\end{lemma}
{\Proof}When $n\equiv1\pmod{6}$, there is an $(n,\{3,4\},1)$-PDF by Lemma \ref{yiwei34}.

When $n\equiv3\pmod{6}$, there is a $(33,\{3,5^{*}\},1)$-PDF and an  $(n,\{3,4^{*},5^{*}\},1)$-PDF for $n \geq 39$ by Lemma \ref{n-345}.  Besides, Equation (\ref{nece}) has only one solution $(x,y,z)=(0,0,1)$ for $n=21$, but there is no $(21,5,1)$-PDF by Lemma \ref{n,3,1}. Similarly, Equation (\ref{nece}) has only one solution $(x,y,z)=(1,0,1)$ for $n=27$, there is no $(27,\{3,5^{*}\},1)$-PDF by Lemma \ref{n-345}.

When $n\equiv5\pmod{6}$, there is an $(n,\{3,4,5\},1)$-PDF for $n\geq59$ by Lemma \ref{n-345}. For $n=41$, Equation (\ref{nece}) has only one solution $(x,y,z)=(0,0,2)$. There is no $(41,5,1)$-PDF by Lemma \ref{n,3,1}. With the help of computer exhaustive search, there is no $(n,\{3,4,5\},1)$-PDF for $n\in\{47,53\}$.

From the above results and Lemma \ref{necessary}, the conclusion is obtained.
\qed

\subsection{Generalized $(n\times m,\{3,4,5\},1)$-PDFs}

In this subsection, we always denote  $D=\{3,5,9,11,15,17,21,23,27,29,35,41,47,53\}$.

\begin{lemma}\label{5SPGDD}
There exists a $\{3,4\}$-SPGDD of type $m^4$ and a $\{3,4,5\}$-SPGDD of type $m^5$ for $m\equiv1\pmod{2}$ and $m\notin D$.
\end{lemma}
{\Proof} There exist a $\{3,4\}$-MGDD of type $k^4$ and a $\{3,4,5\}$-MGDD of type $k^5$ for $k\in\{3,4,5\}$ by Lemma \ref{MGDD}. Applying Construction \ref{MGDD-SPMGDD} with an $(m,\{3,4,5\},1)$-PDF in Lemma \ref{yiwei}, we get a $\{3,4\}$-SPMGDD of type $m^4$ on $A\times [m]$ and a $\{3,4,5\}$-SPMGDD of type $m^5$ on $B\times [m]$, where $m\equiv1\pmod{2}$, $m\notin D$, $A$ is a set of $4$ points, $B$ is a set of $5$ points, and denote the set of base blocks by $\mathcal{A}$ and $\mathcal{B}$, respectively. Let
$$\mathcal{C}=\mathcal{A}\cup\{A\times\{0\}\},~~~~
 \mathcal{D}=\mathcal{B}\cup\{B\times\{0\}\}.$$
Then $\mathcal{C}$ forms a $\{3,4\}$-SPGDD of type $m^4$ and $\mathcal{D}$ forms a $\{3,4,5\}$-SPGDD of type $m^5$.\qed

\begin{lemma}\label{noD}
There exists a generalized $(n\times m,\{3,4,5\},1)$-PDF for $n,m\equiv1\pmod{2}$ and $n,m\notin D$.
\end{lemma}
{\Proof}For $n,m\equiv1\pmod{2}$ and  $n,m\notin D$, there is an $(n,\{3,4,5\},1)$-PDF and an $(m,\{3,4,5\},1)$-PDF by Lemma \ref{yiwei}. By Lemma \ref{5SPGDD} and Lemma \ref{345PDM}, there is a $\{3,4,5\}$-SPGDD of type $m^k$ for $k\in\{3,4,5\}$. Then we get a generalized $(n\times m,\{3,4,5\},1)$-PDF from  Construction \ref{budong} and Construction \ref{N-NxM-PDP}. \qed

It is clear that the existence of a generalized $(n\times m,K,1)$-PDF implies the existence of a generalized $(m\times n,K,1)$-PDF. In what follows,  we  only need to consider the existence of generalized $(n\times m, \{3,4,5\},1)$-PDFs for $m\equiv1\pmod{2}$ and $n\in D$ by Lemma \ref{noD}.

\begin{lemma}\label{23-29-35}
There exists an $(n,3,1)$-PDP with leave $[h]^r$, where ~$(n,h,r)\in\{(23,5,4), (29,5,$ $3),(35,5,3),(41,11,4),(47,11,4),(53,11,3)\}$.
\end{lemma}
{\Proof}All base blocks of an $(n,3,1)$-PDP with leave $[h]^r$ are listed below.\vspace{0.2cm}

$\begin{array}{llllllll}
n=$23$: &$\{0,1,11\}$, &$\{0,2,7\}$,&$\{0,3,9\}$.&&&&\\

n=$29$: &$\{0,2,14\}$, &$\{0,5,13\}$, &$\{0,7,11\}$, &$\{0,9,10\}$.&&&\\

n=$35$: &$\{0,1,10\}$, &$\{0,2,14\}$, &$\{0,4,17\}$, &$\{0,5,16\}$, &$\{0,8,15\}$.&&\\

n=$41$: &$\{0,1,6\}$, &$\{0,2,15\}$, &$\{0,3,17\}$, &$\{0,7,18\}$,&$\{0,9,19\}$.&&\\

n=$47$: &$\{0,2,23\}$, &$\{0,6,17\}$, &$\{0,7,22\}$, &$\{0,9,14\}$, &$\{0,10,13\}$, &$\{0,18,19\}$.&\\

n=$53$:&$\{0,1,24\}$,&$\{0,2,21\}$,&$\{0,5,25\}$, &$\{0,8,22\}$, &$\{0,10,26\}$, &$\{0,11,18\}$, &$\{0,13,17\}$.
\end{array}$

\qed

\begin{lemma}\label{3xm}
There is no generalized $(3\times m,\{3,4,5\},1)$-PDF for $m\equiv1\pmod{2}$.
\end{lemma}
{\Proof}Let $m=2a+1, a\geq 1$. Suppose that there exists a generalized $(3\times m,\{3,4,5\},1)$-PDF with base block set $\mathcal{B}$. Denote the base block set of size  $i$ in $\mathcal{B}$ by $\mathcal{B}_i$ and $|\mathcal{B}_i|=x_i$ for $i\in \{3,4,5\}$, respectively. Then we have
\begin{equation}\label{cha}
6x_3+12x_4+20x_5=6a+2.
\end{equation}

For each $B\in\mathcal{B}$, let
$\Delta B_0=\{b: (0,b)\in\Delta B, b>0\}$. By  listing  all of the possible forms of the  base blocks, we have   that $|\Delta B_{0}|=1$ or $3$ for $B\in  \mathcal{B}_3$, $|\Delta B_{0}|=2,3$ or $6$ for $B\in \mathcal{B}_4$, $|\Delta B_{0}|=4,6$ or $10$ for $B\in \mathcal{B}_5$.
Denote
$$\Delta\mathcal{B}_0=\bigcup \limits_{B\in \mathcal{B}}\Delta B_{0}.$$

It is not difficult to prove that Equation (\ref{cha}) has a non-negative integer solution $(x_3,x_4,x_5)$ if and only if $x_5=1+3k, k\geq0$. Moreover, if Equation (\ref{cha}) has a solution, then the solution has the form $(x_3,x_4,x_5)=(\delta-2j,j,x_5)$, where $\delta=\frac{6a+2-20x_5}{6}$. Due to $|\Delta B_{0}|\geq1$ for $B\in  \mathcal{B}_3$, $|\Delta B_{0}|\geq2$ for $B\in  \mathcal{B}_4$,~$|\Delta B_{0}|\geq4$ for $B\in  \mathcal{B}_5$,  we have
\begin{align*}\label{3Xm}
  |\Delta\mathcal{B}_0|&\geq(\delta-2j)\cdot1+j\cdot2+x_5\cdot4
  =a+\frac{2x_5+1}{3}>a.
\end{align*}
However, we know $|\Delta\mathcal{B}_0|=a$ by the property of perfect difference families, which  is a contradiction.

Thus, there is no generalized $(3\times m,\{3,4,5\},1)$-PDF for any $m\equiv1\pmod{2}$.\qed

For $n\in\{23,29,35,41,47,53\}$, $m\equiv1\pmod{2}, m>3$, applying Construction \ref{N-NxM-PDP} with  an $(n,3,1)$-PDP in Lemma \ref{23-29-35} and a PDM$(3,m)$ in Lemma \ref{345PDM}, we get a generalized $(n\times m,3,1)$-PDP with leave $[5]^r\times[m]$ or $[11]^r\times[m]$ for some positive integers $r$. Therefore, we first need to consider the existence of a generalized $(n\times m,\{3,4,5\},1)$-PDF for $n\in \{5,9,11,15,17,21,27\}$ and $m\equiv1\pmod{2}, m>3$.

Next, we give some  $(m,3,1)$-PDPs with   special   leaves by using Langford sequences.
\begin{lemma}\label{origin}
There exists an $(m,3,1)$-PDP with leave $[h]$ for any $(m,h)$ satisfying one of the following conditions:\vspace{0.1cm}

$(1)$~$m\equiv1,7\pmod{24}$, and $m\geq49$, $h=7$ or $m\geq175$, $h=25;$

$(2)$~$m\equiv3,21\pmod{24}$, and $m\geq21$, $h=3$ or $m\geq147$, $h=21$ or $m\geq315$, $h=45;$

$(3)$~$m\equiv5,11\pmod{24}$, and $m\geq35$, $h=5$ or $m\geq77$, $h=11;$

$(4)$~$m\equiv9,15\pmod{24}$, and $m\geq63$, $h=9$ or $m\geq105$, $h=15;$

$(5)$~$m\equiv13,19\pmod{24}$, and $m\geq91$, $h=13$ or $m\geq133$, $h=19;$

$(6)$~$m\equiv17,23\pmod{24}$, and $m\geq119$, $h=17$ or $m\geq161$, $h=23$.
\end{lemma}

{\Proof}It is certified by Lemma \ref{LDjieguo} and Construction \ref{LD-PDP}. \qed\vspace{0.2cm}

In what follows, we discuss the existence of generalized $(n\times m,\{3,4,5\},1)$-PDFs according to different values of $n$. For each $n\in \{5,9,11,15,17,21,27\}$, we first need to find the appropriate $(m,3,1)$-PDPs with leave $[h]$ in Lemma \ref{origin}, to ensure the existence of generalized $(n\times h,\{3,4,5\},1)$-PDFs, and consequently to  get the infinite classes of generalized $(n\times m,\{3,4,5\},1)$-PDFs by applying recursive constructions.

\normalsize
\begin{lemma}\label{7,9,13}
There is no generalized $(5\times m,\{3,4,5\},1)$-PDF for $m\in\{7,9\}$.
\end{lemma}
{\Proof}It is known that a generalized $(n\times m,K,1)$-PDP implies an  $(nm,K,1)$-PDP with the same number of base blocks. There is no $(35,\{3,4,5\},1)$-PDF by Lemma \ref{yiwei}, then there is no generalized $(5\times 7,\{3,4,5\},1)$-PDF.

Suppose that there exists a generalized $(5\times 9,\{3,4,5\},1)$-PDF with  base block set $\mathcal{B}$. Denote the set of base blocks of size $i$ in $\mathcal{B}$ by  $\mathcal{B}_i$ and $|\mathcal{B}_i|=x_i$ for $i\in \{3,4,5\}$, respectively. Further, we define
 $$\Delta B_{i,0}=\{(0,y)\in \Delta B: B\in \mathcal{B}_i, y>0\}, ~~~\Delta B_{i,j}=\{(j,y)\in \Delta B: B\in\mathcal{B}_i\},~j=1,2,$$
and
 $$\Delta\mathcal{F}_{j}=\bigcup \limits_{B\in \mathcal{B}}(\Delta B_{i,j}),~j=0,1,2.$$
 From the property of perfect difference families, we have $|\Delta\mathcal{F}_{0}|=4, |\Delta\mathcal{F}_{1}|=|\Delta\mathcal{F}_{2}|=9$.
By listing all of the possible forms of the base blocks, we have \vspace{0.2cm}

$(|\Delta B_{3,0}|,|\Delta B_{3,1}|,|\Delta B_{3,2}|)=(3,0,0)$ or $(1,2,0)$ or $(1,0,2)$ or $(0,2,1);$\vspace{0.1cm}

$(|\Delta B_{4,0}|,|\Delta B_{4,1}|,|\Delta B_{4,2}|)=(6,0,0)$ or $(3,3,0)$ or $(3,0,3)$ or $(2,4,0)$ or $(2,0,4)$ or

\hspace{4.85cm}$(1,3,2)$ or $(1,4,1);$\vspace{0.1cm}

$(|\Delta B_{5,0}|,|\Delta B_{5,1}|,|\Delta B_{5,2}|)=(10,0,0)$ or $(6,4,0)$ or $(6,0,4)$ or $(4,6,0)$ or $(4,0,6)$ or

\hspace{4.85cm}$(3,4,3)$ or $(3,6,1)$ or $(2,6,2)$ or $(2,4,4)$.\vspace{0.2cm}

On the other side, we have $|\Delta\mathcal{B}|=6x_3+12x_4+20x_5=44$.  It is straightforward to check that the equation  has three solutions $(x_3,x_4,x_5)=(4,0,1)$ or $(2,1,1)$ or $(0,2,1)$. $(x_3,x_4,x_5)=(4,0,1)$ is impossible, since $|\Delta B_{3,1}|$ and $|\Delta B_{5,1}|$ are even. If $(x_3,x_4,x_5)=(2,1,1)$, since $|\Delta\mathcal{F}_{1}|$ is odd and $|\Delta B_{5,0}|\geq2$, we have $(|\Delta B_{4,0}|,$ $|\Delta B_{4,1}|,|\Delta B_{4,2}|)=(1,3,2)$. Thus, we have $(|\Delta B_{5,0}|,|\Delta B_{5,1}|,$ $|\Delta B_{5,2}|)=(3,4,3)$ or $(3,6,1)$ or $(2,6,2)$ or $(2,4,4)$, which are also impossible. If $(x_3,x_4,x_5)=(0,2,1)$, note that $|\Delta B_{4,0}|\geq1$ and $|\Delta B_{5,0}|\geq2$, then we have $(|\Delta B_{4,0}|,$ $|\Delta B_{4,1}|,|\Delta B_{4,2}|)=(1,3,2)$ or $(1,4,1)$, $(|\Delta B_{5,0}|,|\Delta B_{5,1}|,$ $|\Delta B_{5,2}|)=(2,6,2)$ or $(2,4,4)$, which are contradictions.
\qed

\begin{lemma}\label{5xm}
There exists a generalized $(5\times m,\{3,4,5\},1)$-PDF for any $m\equiv1\pmod{2}$, $m>3$ and $m\notin\{7,9,13,45\}$.
\end{lemma}
{\Proof}When $m\equiv5\pmod{6}$, we have $5m\equiv1\pmod{6}$, then there is a generalized $(5\times m, \{3,4\},1)$-PDF by Theorem \ref{result34}. Next, we consider the case $m\equiv1,3\pmod{6}$. We classify   the values $m>3$ by modulo 24 and show them in Table 1. For convenience, the input generalized $(5\times h,\{3,4,5\},1)$-PDFs are abbreviated as $5\times h$ listed in the last column of  Table 1.

\renewcommand\arraystretch{1.3}%
\begin{table}[H]\caption{Generalized $(5\times m,\{3,4,5\},1)$-PDFs}\label{biao1} \centering
$\scalebox{0.7}{
\begin{tabular}{|l|l|l|l|}
\hline
    $m\equiv\pmod{24}$     &the values of $m$  & leave $[h]^r$ & input designs  \\ \hline
\multicolumn{1}{|c|}{} & $m\geq175$ & $[25]$  & $5\times25$  \\ \cline{2-3}
                       &$m\in\{79,97,103,121,127,145,151,169\}$  & $[25]^3$ &   \\ \cline{2-4}
~~~~~~~$1,7$           & $m=73$ & $[19]^4$ & $5\times19$  \\ \cline{2-4}
                       & $m\in\{31,49,55\}$ & $[5]^4$ & $5\times5$   \\ \cline{2-4}
                       & $m=\{7,25\}$ & $[m]$ & $5\times m$   \\ \hline
                       & $m\geq147$ & $[21]$ &    \\ \cline{2-3}
                       & $m\in\{69,75,93,99,117,123,141\}$ & $[21]^3$ & $5\times21$   \\ \cline{2-4}
~~~~~~$3,21$           & $m=51$ & $[15]^3$ & $5\times15$   \\ \cline{2-4}
                       & $m\in\{21,27,45\}$ & $[m]$ & $5\times m$  \\ \hline
                       & $m\geq105$ & $[15]$ & $5\times15$   \\ \cline{2-3}
~~~~~~$9,15$           & $m\in\{57,63,81,87\}$ & $[15]^3$  &\\ \cline{2-4}
                       & $m\in\{9,15,33,39\}$ & $[m]$ & $5\times m$   \\ \hline
                       & $m\geq133$ & $[19]$ & $5\times19$   \\ \cline{2-3}
                       & $m\in\{61,67,85,91,109,115\}$ & $[19]^3$ & \\ \cline{2-4}
~~~~~~$13,19$   & $m\in\{37,43\}$ & $[5]^5$ & $5\times 5$   \\ \cline{2-4}
                       & $m=\{13,19\}$ & $[m]$ & $5\times m$   \\ \hline
\end{tabular}}$
\end{table}

For $m$ in Table 1, there is an $(m,\{3,4,5\},1)$-PDP with leave $[h]^r$ by Lemma \ref{origin} and Appendix \ref{app-C}, \ref{app-D}. By Lemma  \ref{345PDM}, there is a PDM$(k,5)$, i.e., a $k$-SPGDD of type $5^k$ for $k\in\{3,4,5\}$, then we get a generalized $(5\times m,\{3,4,5\},1)$-PDP with leave $[5]\times [h]^r$ by Construction \ref{N-NxM-PDP}. Applying Construction \ref{budong}, we only need to construct the input   generalized $(5\times h,\{3,4,5\},1)$-PDFs, where $h\in\{5,7,9,13,15,19,21,25,27,33,39,45\}$.
The generalized   $(5\times5,\{3,4\},1)$-PDF  has been given in Theorem \ref{result34}, and there is no generalized $(5\times m,\{3,4,5\},1)$-PDFs for $m\in\{7,9\}$ by Lemma \ref{7,9,13}. For $h\in\{15,19,21,25,27,33,39\}$, the existence of generalized $(5\times h,\{3,4,5\},1)$-PDFs see Appendix A. \qed

\begin{lemma}\label{9xm}
There is a generalized $(n\times m,\{3,4,5\},1)$-PDF for $n\in\{9,11\}, m\equiv1\pmod{2}$,  $m>3$ and $(n,m)\notin\{(9,5),(9,35),(11,19),(11,27)\}$.
\end{lemma}

{\Proof}Similar to Table 1, we construct Table 2 below.

For $m$ in Table 2, there is an $(m,\{3,4\},1)$-PDP with leave $[h]^r$ from Lemma \ref{23-29-35}, Lemma \ref{origin} and Appendix \ref{app-C}. For $n\in\{9,11\}$,  there is a $3$-SPGDD of type $n^3$ by Lemma \ref{345PDM}  and a $\{3,4\}$-SPGDD of type $n^4$ by Example \ref{9SPGDD} and Example \ref{11SPGDD}. Then we get a generalized $(n\times m,\{3,4\},1)$-PDP with leave $[n]\times [h]^r$ by Construction \ref{N-NxM-PDP}. Applying  Construction \ref{budong}, we only need to construct the input generalized $(n\times h,\{3,4,5\},1)$-PDFs,  where $n\in\{9,11\}, h\in\{5,7,9,11,13,15,17,19,21,23,27,29,35\}$.

\renewcommand\arraystretch{1.3}%
\begin{table}[H]\caption{Generalized $(n\times m,\{3,4,5\},1)$-PDFs for $n\in\{9,11\}$}\label{biao2} \centering
$\scalebox{0.7}{
\begin{tabular}{|l|l|l|l|}
\hline
    $m\equiv\pmod{24}$     &the values of $m$  & leave $[h]^r$  & input designs  \\ \hline
\multicolumn{1}{|c|}{} & $m\geq49$ & $[7]$  &   \\ \cline{2-3}
 ~~~~~~~$1,7$          &$m\in\{25,31\}$  & $[7]^3$ &  $n\times7$  \\ \cline{2-3}
                       & $m=7$ & $[7]$ &  \\ \hline
                       & $m\geq147$ & $[21]$ &  $n\times21$ \\ \cline{2-4}
~~~~~~$3,21$           & $m\in\{45,51,69,75,93,99,117,123,141\}$ & $[9]^3$ & $n\times9$  \\ \cline{2-4}
                       & $m\in\{21,27\}$ & $[m]$ & $n\times m$  \\ \hline
                       & $m\geq63$ & $[9]$ & $n\times9$  \\ \cline{2-3}
~~~~~~$9,15$           & $m\in\{33,39,57\}$ & $[9]^3$ &    \\ \cline{2-4}
                       & $m=\{9,15\}$ & $[m]$ & $n\times m$  \\ \hline
                       & $m\geq91$ & $[13]$ & $n\times13$   \\ \cline{2-4}
~~~~~~$13,19$          & $m\in\{37,43,61,67,85\}$ & $[7]^3$ & $n\times7$  \\ \cline{2-4}
                       & $m=\{13,19\}$ & $[m]$ & $n\times m$ \\ \hline
                       & $m\geq77$ & $[11]$ & $n\times11$  \\ \cline{2-3}
~~~~~~$5,11$           & $m\in\{53,59\}$ & $[11]^3$ &    \\ \cline{2-4}
                       & $m\in\{5,11,29,35\}$ & $[m]$ & $n\times m$ \\ \hline
                       & $m\geq119$ & $[17]$ & $n\times17$  \\ \cline{2-4}
~~~~~~$17,23$          & $m\in\{65,71,89,95,113\}$ & $[11]^3$ & $n\times11$  \\ \cline{2-3}
                       & $m\in\{41,47\}$ & $[11]^4$ &    \\ \cline{2-4}
                       & $m\in\{17,23\}$ & $[m]$ & $n\times m$  \\ \hline
\end{tabular}}$
\end{table}

For $n=9, h\in\{7,9,11,13,15,17,19,21,23,29\}$, the  generalized $(n\times h,\{3,4,5\},1)$-PDFs see Appendix A. There is no generalized $(9\times 5,\{3,4,5\},1)$-PDF by Lemma \ref{7,9,13}.
We here construct a generalized $(9\times27,\{3,4,5\},1)$-PDF, whose base blocks are divided into two parts. For the first part, we know  $\{\{0,1,3,11\}\}$ is a $(27,4,1)$-PDP with leave $N=\{0,\pm4,\pm5,\pm6,\pm7,\pm9,\pm12,\pm13\}$. There is a $\{3,4\}$-SPGDD of type $9^4$ by Example \ref{9SPGDD}. Applying Construction \ref{N-NxM-PDP}, we get a generalized $(9\times27,\{3,4\},1)$-PDP with leave $[9]\times N$. The second part consists of the base blocks of a generalized $([9]\times N,\{3,4,5\},1)$-PDF, which are listed below.\vspace{-0.4cm}

{\small$$\begin{array}{llllllll}
&$\{(0, 0), (0, 7), (4, 13)\}$,~~~~
&$\{(0, 9), (1, 4), (4, 0)\}$,~~
&$\{(0, 5), (3, 12), (4, 0)\}$,\\
&$\{(0, 6), (1, 12), (3, 0)\}$,~~~~
&$\{(0, 12), (3, 7), (4, 0)\}$,~~
&$\{(0, 0), (2, 5), (0, 12)\}$,\\
&$\{(0, 13), (1, 7), (4, 0)\}$,~~~~
&$\{(0, 4), (2, 0), (4, 9)\}$,~~
&$\{(0, 0), (2, 4), (2, 13)\}$,\\
&$\{(0, 7), (1, 12), (4, 0)\}$,~~~~
&$\{(0, 0), (3, 5), (4, 9)\}$,~~
&$\{(0, 0), (1, 0), (3, 6), (3, 12)\}$,\\
&$\{(0, 13), (1, 0), (4, 9)\}$,~~~~
&$\{(0, 0), (3, 0), (4, 7), (4, 12)\}$,~~
&$\{(0,~6), (2,~0), (2,~13), (4,~0)\}$,
\vspace{-0.2cm}
\end{array}$$
\begin{flushleft}
$\hspace{0.77cm}\{(0,~0), (3,~13), (4,~4)\},\hspace{0.01cm}$
$\hspace{0.7cm}\{(0, 13), (1, 9), (3, 0), (3, 4), (4, 13)\}$.
\end{flushleft}}

For $n=11, h\in\{5,11,17,23,29\}$,  we have $11h\equiv1\pmod{6}$, then there is a generalized $(11\times h, \{3,4\},1)$-PDF by Theorem \ref{result34}. For $h\in\{7,9,13,15,21\}$, the results of PDFs see Appendix A.\qed

\begin{lemma}\label{15xm}
There exists a generalized $(n\times m,\{3,4,5\},1)$-PDF for $n\in\{15,17,21,27\},  m\equiv1\pmod{2}, m>3$ and $(n,m)\notin\{(15,17),(15,21),(15,27),(17,13),(17,15),(17,21),$ $(17,27),$ $(21,13),(21,15),
(21,17),(21,21),(21,27),(27,11),(27,13),(27,15),(27,17),(27,21),(27,27)\}$.
\end{lemma}

{\Proof}Similar to Table 1, we construct Table 3 below.

For $m$ in Table \ref{biao3}, there is an $(m,\{3,4\},1)$-PDP with leave $[h]^r$ by Lemma \ref{23-29-35}, Lemma \ref{origin} and Appendix \ref{app-C}. For $n\in\{15,17,21,27\}$, there is a  $k$-SPGDD of type $n^k$ for $k\in\{3,4\}$ (i.e., a PDM$(k,n)$ by Lemma \ref{345PDM}), then we get a generalized $(n\times m,\{3,4\},1)$-PDP with leave $[n]\times [h]^r$ by Construction \ref{N-NxM-PDP}. Applying Construction \ref{budong}, we only need to construct the input generalized $(n\times h,\{3,4,5\},1)$-PDFs, where $n\in\{15,17,21,27\}$, $h\in\{5,7,9,11,13,15,17,19,21,23,27,45\}$.

\renewcommand\arraystretch{1.3}%
\begin{table}[H]\caption{Generalized $(n\times m,\{3,4,5\},1)$-PDFs for $n\in\{15,17,21,27\}$}\label{biao3} \centering
$\scalebox{0.7}{
\begin{tabular}{|l|l|l|l|}
\hline
    $m\equiv\pmod{24}$     &the values of $m$ & leave $[h]^r$ & input designs  \\ \hline
\multicolumn{1}{|c|}{} & $m\geq49$ & $[7]$  &  \\ \cline{2-3}
 ~~~~~~~$1,7$          &$m\in\{25,31\}$  & $[7]^3$ &  $n\times7$   \\ \cline{2-3}
                       & $m=7$ & $[7]$ &  \\ \hline
                       & $m\geq315$ & $[45]$ & $n\times45$   \\ \cline{2-4}
                       & $m\in\{45,51,69,75,93,99,117,123,141,147,165,171,$ & $[9]^3$ & $n\times9$ \\
~~~~~~$3,21$           & ~~~~~~~~$189,195,213,219,237,243,261,267,285,291,309\}$ &  &   \\ \cline{2-4}
                       & $m\in\{21,27\}$ & $[m]$ & $n\times m$  \\ \hline
                       & $m\geq63$ & $[9]$ & $n\times9$  \\ \cline{2-3}
~~~~~~$9,15$           & $m\in\{33,39,57\}$ & $[9]^3$ &    \\ \cline{2-4}
                       & $m\in\{9,15\}$ & $[m]$ & $n\times m$ \\ \hline
                       & $m\geq133$ & $[19]$ & $n\times19$  \\ \cline{2-4}
~~~~~~$13,19$          & $m\in\{37,43,61,67,85,91,109,115\}$ & $[7]^3$ & $n\times7$ \\  \cline{2-4}
                       & $m\in\{13,19\}$ & $[m]$ & $n\times m$  \\ \hline
                       & $m\geq35$ & $[5]$ & $n\times5$  \\ \cline{2-3}
~~~~~~$5,11$           & $m=29$ & $[5]^3$ &    \\ \cline{2-4}
                       & $m\in\{5,11\}$ & $[m]$ & $n\times m$  \\ \hline
                       & $m\geq161$ & $[23]$ & $n\times23$  \\ \cline{2-4}
                       & $m\in\{41,47,65,71,89,95,113,119,137,143\}$ & $[5]^3$ & $n\times5$  \\ \cline{2-3}
~~~~~~$17,23$          & $m=23$ & $[5]^4$ &   \\ \cline{2-4}
                       & $m=17$ & $[17]$ & $n\times 17$ \\ \hline
\end{tabular}}$
\end{table}
 For $n\in\{15,17,21,27\}$, when $h=19$, we here construct a generalized $(n\times19,\{3,4,5\},1)$-PDF, whose base blocks are divided into two parts. For the first part, we know  $\{\{0,3,5,9\}\}$ is a $(19,4,1)$-PDP with leave $M=\{0,\pm1,\pm7,\pm8\}$. Applying Construction \ref{N-NxM-PDP} with a $4$-SPGDD of type $n^4$ in Lemma \ref{345PDM}, we get a generalized $(n\times19,4,1)$-PDP with leave $[n]\times M$. The second part consists of the base blocks of a generalized $([n]\times M,\{3,4,5\},1)$-PDF, see   Appendix \ref{app-B}. When $h\in\{23,45\}$,  we only  need to construct an input  generalized $(n\times5,\{3,4,5\},1)$-PDF and a generalized $(n\times9,\{3,4,5\},1)$-PDF from Table 3.

To be specific,  for $n=17$, $h\in\{5,11,17,23,29\}$,  we have $17h\equiv1\pmod{6}$,  then there is also a generalized $(17\times h, \{3,4\},1)$-PDF by Theorem \ref{result34}. The generalized $(27\times 9, \{3,4,5\},1)$-PDF exists  by Lemma \ref{9xm}.
For $n=15$ and $ h\in\{5,7,9,11,13,15\}$, $n=17$ and $h\in\{7,9\}$, $n=21$ and $ h\in\{5,7,9,11\}$, $n=27$ and  $h\in\{5,7\}$, the results of PDFs see Appendix A.\qed

Based on the recursive constructions and the non-existence of some generalized $(n\times m,\{3,4,5\},1)$-PDFs for $n\in \{5,11\}$, we also need to consider the existence of generalized $(n\times m,\{3,4,5\},1)$-PDFs for $n\in \{23,29,35,41,47,53\}$.

\begin{lemma}\label{partD2}
For $n\in \{23,29,35,41,47,53\}$, there exists a generalized $(n\times m,\{3,4,5\},1)$-PDF for $m\equiv1\pmod{2}$ and $m>3$, except for $n=35$ and $m\in\{7,9,13,45\}$, or $n\in\{23,29\}$ and $m\in\{7,13,45\}$.
\end{lemma}
{\Proof}When $n\in\{23,29,35\}$, $m\equiv1\pmod{2}$ and $m>3$, there is an  $(n,3,1)$-PDP with leave $[5]^r$ for some positive integers $r$
by Lemma \ref{23-29-35}. There is a 3-SPGDD of type $m^3$ by Lemma \ref{345PDM}. Then we get a generalized $(n\times m,3,1)$-PDP with leave $[5]^r\times m$ by Construction \ref{N-NxM-PDP}. For $m\notin\{7,9,13,45\}$, there is a generalized $(5\times m,\{3,4,5\},1)$-PDF   from  Lemma \ref{5xm}, then we get a generalized $(n\times m,3,1)$-PDF by Construction \ref{budong}. Further, the existence of a generalized $(n\times9,\{3,4,5\},1)$-PDF for $n\in\{23,29\}$ see Appendix A.

When $n\in\{41,47,53\}$, $m\equiv1\pmod{2}$ and $m>3$, there is an $(n,3,1)$-PDP with leave $[11]^r$ for some positive integers $r$ by Lemma \ref{23-29-35} and a 3-SPGDD of type $m^3$ by Lemma \ref{345PDM}. Then we get a generalized $(n\times m,3,1)$-PDP with leave $[11]^r\times m$ by Construction \ref{N-NxM-PDP}.
For $m\notin\{19,27\}$, there is a generalized $(11\times m,\{3,4,5\},1)$-PDF from Lemma \ref{9xm}. For $m\in\{19,27\}$, applying Construction \ref{N-NxM-PDP} with an $(n,\{3,4\},1)$-PDP with leave $[5]^3$ in Appendix \ref{app-C} and a $k$-SPGDD of type $m^k$ for $k\in\{3,4\}$ in Lemma \ref{345PDM}, we get a generalized $(n\times m,3,1)$-PDP with leave $[5]^r\times m$. A generalized $(5\times m,\{3,4,5\},1)$-PDF exists for $m\in\{19,27\}$ by Lemma \ref{5xm}. Applying Construction \ref{budong}, we get all the generalized PDFs in this situation. \qed

\emph{Proof of Theorem \ref{result345}}

Combining the results of Lemmas \ref{noD}, \ref{3xm}, \ref{7,9,13}--\ref{partD2}, we complete the proof.

\section{Applications}

In this section, we discuss the application of generalized perfect difference families to variable-weight geometric orthogonal codes.
Geometric orthogonal codes were first introduced by Doty and Winslow \cite{Doty}, who  used geometric orthogonal codes to design sets of macrobonds in $3$D DNA origami so as to reduce undesirable bonding arising from misalignment and mismatches.
Rothemund \cite{Rothemund} indicated that when some external conditions,  such as temperature, salinity, density and concentration are the same, compared with a single bonds strength, different bonds strengths (codewords weights) improve the efficiency of constructing nanostructures. Therefore,  it is very important to study variable-weight GOCs.

Let $\mathbb{Z}$ be the set of integers and $n, m, \lambda_a, \lambda_c$ be positive integers. An {\em  $(n\times m,K,\lambda_a,\lambda_c)$-geometric orthogonal code}, briefly  $(n\times m,K,\lambda_a,\lambda_c)$-GOC, is a collection $\cal{C}$ of subsets (called \emph{codewords} or \emph{macrobonds}) of $[0,n-1]\times[0,m-1]$ of sizes from a set of positive integers $K$ such that:
\begin{itemize}
\item[$(1)$] (the aperiodic auto-correlation):
$|B\cap(B+(s,t))|\leq\lambda_a$ for all $B\in\cal{C}$ and every $(s,t)\in \mathbb{Z}\times \mathbb{Z}\setminus\{(0,0)\}$;
\item[$(2)$] (the aperiodic cross-correlation):
$|A\cap(B+(s,t))|\leq\lambda_c$ for all $A,B\in\cal{C}$ with $A\neq B$ and every $(s,t)\in\mathbb{Z}\times \mathbb{Z}$,
\end{itemize}
where $B+(s,t)=\{(x+s,y+t):(x,y)\in B\}$. When $\lambda_a=\lambda_c=\lambda$, the notation of $(n\times m,K,\lambda)$-GOC is employed for short. Such GOCs are called \emph{constant-weight} when $|K|=1$ or \emph{variable-weight} when $|K|>1$.

Let $\mathcal{C}$ be a collection of subsets of $[0,n-1]\times[0,m-1]$ of sizes from $K$. It is not convenient to check
the correctness of Conditions $(1)$ and $(2)$. But fortunately, the difference method is very efficient to describe an $(n\times m,K,\lambda)$-GOC for $\lambda=1$. Let $\Delta B$ be the list of differences from $B$ for any $B\in\mathcal{C}$ and $\Delta\mathcal{C}=\bigcup_{B\in\mathcal{C}}\Delta B$. It is readily checked that $\mathcal{C}$ constitutes an $(n\times m,K,1)$-GOC if
$\Delta\mathcal{C}$ covers every element of $[2n-1]\times[2m-1]$ at most once.
Furthermore, $\mathcal{C}$ is called a {\em perfect } $(n\times m,K,1)$-GOC if $\Delta\mathcal{C}=[2n-1]\times[2m-1]\backslash\{(0,0)\}$.
\newpage

\begin{lemma}\label{GOC-PDF}
A perfect $(n\times m,K,1)$-GOC is equivalent to a generalized $((2n-1)\times(2m-1),K,1)$-PDF.
\end{lemma}
{\Proof}Clearly, a perfect $(n\times m,K,1)$-GOC is a generalized $((2n-1)\times(2m-1),K,1)$-PDF.
Conversely, suppose $\cal{B}$ is a generalized $((2n-1)\times(2m-1),K,1)$-PDF. For each $B=\{(x_i,y_i):1\leq i\leq k\}\in\cal{B}$, let $C_B=\{(x_i-x_0,y_i-y_0):1\leq i\leq k\}$, where $x_0=\min\{x_i:1\leq i\leq k\}$ and $y_0=\min\{y_i:1\leq i\leq k\}$. Let  ${\cal{C}}=\bigcup_{B\in\cal{B}}C_B$. Then $\cal{C}$ is a perfect $(n\times m,K,1)$-GOC. \qed

By Theorems \ref{result34}, \ref{result345} and Lemma \ref{GOC-PDF}, we get some results of perfect GOCs as follows.

\begin{theorem}\label{GOC34-result1}
There exists a variable-weight perfect $(n\times m,\{3,4\},1)$-GOC if and only if  $n,m\equiv0\pmod{3}$ or $n,m\equiv1\pmod{3}$.
\end{theorem}

\begin{theorem}\label{GOC345-result2}
 There exists a variable-weight perfect $(n\times m,\{3,4,5\},$ $1)$-GOC if and only if $n,m\neq2$ and $\{n,m\}\notin\{\{1,d\},\{3,4\},\{3,5\}\}$, where $d\in\{3,5,6,8,9,11,$ $12,14,15,18,21,24,$ $27\}$, and possibly except for
$\{n,m\}\in\{\{3,7\},\{3,23\},\{4,12\},\{4, 15\},\{4,$ $18\},\{5,18\},\{6,10\},
$ $\{6,14\},\{7,9\},~\{7,11\},~\{7,12\},\{7,14\},\{7,15\},\{7,18\},\{8,9\}, \{8, 11\},\{8,$ $14\},\{9,11\},\{9,14\},$ $\{11,11\},
\{11,14\}, \{12,23\},\{14,14\},\{15,23\},\{18,23\}\}$.
\end{theorem}

\section{Concluding remarks}

Based on the equivalence between a variable-weight perfect GOC and a generalized PDF, generalized $(n\times m,K,1)$-PDFs were investigated in this paper. By using some auxiliary designs and several recursive constructions, the existence of generalized $(n\times m,\{3,4\},1)$-PDFs is completely solved. The existence of a generalized $(n\times m, \{3,4,5\}, 1)$-PDF is also completely solved possibly except for a few values. As a consequence, the existence of corresponding variable-weight perfect $(n\times m,\{3,4\},1)$-GOCs and $(n\times m,\{3,4,5\},1)$-GOCs are obtained. Some new constructions are expected in the further research to completely resolve the existence problem of generalized $(n\times m, \{3,4,5\}, 1)$-PDFs and study the other cases of generalized $(n\times m,K,1)$-PDFs such as for $K=\{3,5\}$.

For variable-weight $(n\times m,K,1)$-GOC, ${\cal C}$, we can also  limit that: there are exactly $q_i |{\cal C}|$ codewords of weight $k_i$, i.e., $q_i$ indicates
the fraction of codewords of weight $k_i$, where $K=\{k_1,k_2,\dots,k_s\},  Q=\{q_1,q_2,\dots,q_s\},$ and $\sum _{i=1}^s q_i=1.$ Actually, this kind of GOCs
were proposed in \cite{Doty} as an open question about variable-weight GOCs.
According to actual needs, studying this kind of GOCs  will be an  interesting and important work.

\newpage
{\centering\title{\bf\Large{ Appendix}}}
\appendix
\section{Generalized $(n\times m,\{3,4,5\},1)$-PDFs}\label{app-A}
The base blocks of a generalized $(n\times m,\{3,4,5\},1)$-PDF are listed below:\\

\noindent1.~$(n,m)=(5,15):$
{\scriptsize \begin{quote}
$\{(0, 1), (1, 0), (2, 7)\},~~
\{(0, 0), (1, 1), (2, 5)\},~~
\{(0, 0), (0, 2), (2, 4)\},~~
\{(0, 5), (1, 2), (2, 0)\},~~$

$\{(0, 6), (1, 0), (2, 0)\},~~
\{(0, 3), (0, 7), (2, 0)\},~~
\{(0, 2), (2, 0), (2, 5)\},~~
\{(0, 1), (0, 4), (1, 7), (2, 0)\},~~$

$\{(0, 0), (0, 6), (0, 7), (1, 2), (2, 7)\}.$
\end{quote}}

\noindent2.~$(n,m)=(5,19):$~
{\scriptsize \begin{quote}
$\{(0, 0), (1, 9), (2, 0)\},~~
\{(0, 1), (1, 0), (2, 7)\},~~
\{(0, 3), (1, 7), (2, 0)\},~~
\{(0, 8), (2, 0), (2, 1)\},$

$\{(0, 9), (1, 4), (2, 0)\},~~
\{(0, 0), (1, 0), (2, 8)\},~~
\{(0, 1), (2, 0), (2, 6)\},~~
\{(0, 0), (0, 4), (1, 2), (2, 7)\},$

$\{(0, 8), (1, 0), (1, 9), (2, 3), (2, 6)\},~~
\{(0, 0), (0, 5), (0, 7), (2, 1), (2, 9)\}.$
\end{quote}}

\noindent3.~$(n,m)=(5,21):$~
{\scriptsize \begin{quote}
$\{(0, 0), (1, 0), (1, 10)\},\hspace{0.2cm}
\{(0, 0), (2, 0), (2, 9)\},\hspace{0.2cm}
\{(0, 0), (2, 1), (2, 8)\},\hspace{0.3cm}
\{(0, 1), (2, 0), (2, 8)\},$

$\{(0, 3), (2, 0), (2, 6)\},\hspace{0.33cm}
\{(0, 0), (2, 2), (2, 6)\},\hspace{0.2cm}
\{(0, 0), (1, 9), (2, 4)\},\hspace{0.3cm}
\{(0, 2), (1, 10), (2, 0)\},$

$\{(0, 4), (1, 9), (2, 0)\},\hspace{0.33cm}
\{(0, 6), (1, 2), (2, 0)\},\hspace{0.2cm}
\{(0, 10), (1, 3), (2, 0)\},\hspace{0.17cm}
\{(0, 8), (2, 0), (2, 3)\},$

$\{(0, 7), (0, 9), (1, 1), (2, 0)\},\hspace{0.2cm}
\{(0, 0), (0, 5), (1, 6), (1, 7), (2, 10)\}.$
\end{quote}}

\noindent4.~$(n,m)=(5,25):$~
{\scriptsize \begin{quote}
$\{(0, 7), (1, 0), (2, 0)\},~~
\{(0, 1), (1, 12), (2, 0)\},~~
\{(0, 10), (1, 1), (2, 0)\},\hspace{0.34cm}
\{(0, 0), (0, 7), (2, 2)\},~~$

$\{(0, 0), (2, 4), (2, 9)\},~~
\{(0, 10), (1, 0), (2, 6)\},~~
\{(0, 11), (1, 0), (2, 12)\},~~
\{(0, 3), (2, 0), (2, 11)\},~~$

$\{(0, 2), (1, 0), (2, 9)\},~~
\{(0, 0), (1, 3), (2, 11)\},~~
\{(0, 0), (0, 8), (1, 10)\},\hspace{0.35cm}
\{(0, 0), (1, 1), (2, 5)\},~~~$

$\{(0, 0), (1, 7), (2, 3), (2, 12)\},~
\{(0, 0), (0, 12), (2, 0), (2, 6), (2, 10)\},~~
\{(0, 8), (0, 9), (0, 11), (1, 3), (2, 0)\}.$
\end{quote}}

\noindent5.~$(n,m)=(5,27):$~
{\scriptsize \begin{quote}
$\{(0, 9), (1, 0), (2, 0)\},~~
\{(0, 0), (1, 13), (2, 0)\},~~
\{(0, 1), (0, 12), (2, 0)\},~~
\{(0, 0), (0, 12), (2, 1)\},~~$

$\{(0, 7), (2, 0), (2, 1)\},~~
\{(0, 2), (1, 11), (2, 0)\},~~
\{(0, 2), (1, 0), (2, 12)\},~~
\{(0, 0), (2, 5), (2, 13)\},~~$

$\{(0, 0), (2, 2), (2, 9)\},~~
\{(0, 10), (1, 0), (2, 6)\},~~
\{(0, 0), (2, 7), (2, 11)\},~~
\{(0, 0), (1, 8), (2, 4)\},~~$

$\{(0, 8), (1, 0), (2, 5)\},~~
\{(0, 5), (2, 0), (2, 13)\},~~
\{(0, 0), (1, 10), (2, 3)\},~~
\{(0, 0), (0, 9), (1, 3)\},~~$

$\{(0, 3), (1, 0), (1, 10)\},~
\{(0, 0), (0, 6), (1, 1), (2, 12)\},~~
\{(0, 8), (0, 10), (0, 13), (1, 12), (2, 0)\}.$
\end{quote}}

\noindent6.~$(n,m)=(5,33):$~
{\scriptsize \begin{quote}
$\{(0, 16), (1, 0), (2, 0)\},~~
\{(0, 0), (0, 15), (2, 0)\},~~
\{(0, 1), (0, 14), (2, 0)\},~~
\{(0, 0), (0, 1), (2, 16)\},~~$

$\{(0, 12), (1, 1), (2, 0)\},~~
\{(0, 0), (0, 11), (2, 1)\},~~
\{(0, 14), (1, 0), (2, 1)\},~~
\{(0, 2), (0, 8), (2, 0)\},~~$

$\{(0, 11), (1, 2), (2, 0)\},~~
\{(0, 3), (1, 15), (2, 0)\},~~
\{(0, 9), (2, 0), (2, 2)\},\hspace{0.33cm}
\{(0, 8), (1, 0), (2, 2)\},~~$

$\{(0, 3), (1, 0), (1, 14)\},~~
\{(0, 0), (0, 16), (1, 3)\},~~
\{(0, 4), (2, 0), (2, 9)\},\hspace{0.33cm}
\{(0, 7), (1, 0), (2, 14)\},~~$

$\{(0, 0), (1, 16), (2, 6)\},~~
\{(0, 4), (1, 0), (1, 12)\},~~
\{(0, 0), (0, 5), (2, 13)\},~~
\{(0, 0), (1, 5), (2, 11)\},~~$

$\{(0, 0), (0, 3), (2, 12)\},~~
\{(0, 0), (0, 4), (1, 13)\},~~
\{(0, 0), (0, 10), (1, 4), (2, 14)\},~~$

$\{(0, 0), (0, 8), (1, 15), (2, 3), (2, 10)\}.$
\end{quote}}

\noindent7.~$(n,m)=(5,39):$~
{\scriptsize \begin{quote}
$\{(0, 1), (0, 16), (2, 0)\},~~
\{(0, 14), (1, 0), (2, 0)\},~~
\{(0, 12), (1, 0), (2, 1)\},~~
\{(0, 0), (1, 2), (2, 16)\},~~$

$\{(0, 0), (2, 1), (2, 17)\},~~
\{(0, 0), (0, 12), (2, 3)\},~~
\{(0, 5), (2, 0), (2, 19)\},~~
\{(0, 0), (0, 4), (2, 12)\},~~$

$\{(0, 3), (1, 13), (2, 0)\},~~
\{(0, 8), (1, 15), (2, 0)\},~~
\{(0, 0), (1, 17), (2, 7)\},~~
\{(0, 2), (1, 6), (2, 0)\},~~$

$\{(0, 8), (1, 0), (2, 19)\},~~
\{(0, 4), (1, 19), (2, 0)\},~~
\{(0, 6), (1, 17), (2, 0)\},~~
\{(0, 0), (0, 10), (1, 3), (2, 19)\},~~$
\end{quote}
\begin{quote}
$\{(0, 0), (0, 18), (2, 0), (2, 5), (2, 6)\},~~~~
\{(0, 0), (0, 3), (0, 14), (1, 12), (2, 18)\},~~
\{(0, 5), (1, 0), (1, 13), (2, 18)\},~~$

$\{(0, 19), (1, 18), (2, 0), (2, 2), (2, 9)\},~~
\{(0, 7), (0, 15), (1, 4), (2, 0), (2, 17)\}.$
\end{quote}}

\noindent8.~$(n,m)=(7,9):$~
{\scriptsize \begin{quote}
$\{(0, 0), (1, 0), (3, 4)\},~~
\{(0, 4), (1, 0), (3, 0)\},~~
\{(0, 0), (0, 3), (3, 0)\},~~
\{(0, 1), (1, 4), (3, 0)\},~~$

$\{(0, 3), (2, 0), (2, 1)\},~~
\{(0, 2), (1, 0), (2, 4)\},~~
\{(0, 0), (3, 1), (3, 3)\},~~
\{(0, 2), (1, 1), (2, 3), (3, 0), (3, 4)\}.$
\end{quote}}

\noindent9.~$(n,m)=(7,11):$~
{\scriptsize \begin{quote}
$\{(0, 0), (1, 0), (3, 5)\},~~
\{(0, 0), (1, 4), (2, 3), (3, 0)\},~~
\{(0, 2), (1, 0), (3, 4)\},~~
\{(0, 5), (1, 0), (3, 2)\},~~$

$\{(0, 1), (0, 5), (3, 0)\},~~
\{(0, 0), (0, 3), (0, 5), (2, 0), (3, 1)\},~
\{(0, 1), (0, 2), (1, 4), (2, 0), (3, 5)\}.$
\end{quote}}

\noindent10.~$(n,m)=(7,15):$~
{\scriptsize \begin{quote}
$\{(0, 0), (1, 0), (3, 7)\},~~
\{(0, 7), (1, 0), (3, 0)\},~~
\{(0, 0), (0, 6), (3, 0)\},~~
\{(0, 1), (0, 5), (3, 0)\},~~$

$\{(0, 7), (1, 1), (2, 0)\},~~
\{(0, 4), (1, 1), (3, 0)\},~~
\{(0, 4), (2, 0), (3, 1)\},~~
\{(0, 2), (1, 6), (3, 0)\},~~$

$\{(0, 0), (0, 7), (2, 2)\},~~
\{(0, 0), (1, 5), (3, 2)\},~~
\{(0, 0), (0, 2), (2, 6)\},~~
\{(0, 0), (1, 6), (2, 1)\},~~$

$\{(0, 0), (0, 3), (3, 6)\},~~
\{(0, 0), (0, 1), (3, 5)\},~~
\{(0, 2), (1, 0), (1, 5), (2, 7), (3, 3)\}.$
\end{quote}}

\noindent11.~$(n,m)=(7,17):$~
{\scriptsize \begin{quote}
$\{(0, 0), (1, 0), (3, 8)\},~~
\{(0, 8), (1, 0), (3, 0)\},~~
\{(0, 0), (0, 7), (3, 0)\},~~
\{(0, 1), (3, 0), (3, 8)\},~~$

$\{(0, 6), (1, 1), (3, 0)\},~~
\{(0, 5), (2, 1), (3, 0)\},~~
\{(0, 2), (0, 8), (2, 0)\},~~
\{(0, 6), (2, 0), (3, 2)\},~~$

$\{(0, 3), (1, 7), (3, 0)\},~~
\{(0, 0), (2, 1), (3, 6)\},~~
\{(0, 0), (2, 3), (2, 7)\},~~
\{(0, 2), (1, 0), (3, 6)\},~~$

$\{(0, 0), (1, 7), (3, 2)\},~~
\{(0, 0), (0, 2), (1, 3), (1, 8), (3, 5)\},~
\{(0, 4), (0, 7), (1, 0), (1, 1), (3, 5)\}.$
\end{quote}}

\noindent12.~$(n,m)=(7,21):$~
{\scriptsize \begin{quote}
$\{(0, 6), (0, 8), (3, 0)\},~~~
\{(0, 0), (0, 7), (3, 2)\},~~
\{(0, 0), (1, 5), (3, 9)\},~~
\{(0, 6), (0, 10), (2, 0)\},~~$

$\{(0, 0), (1, 6), (3, 8)\},~~~
\{(0, 9), (1, 0), (3, 9)\},~~
\{(0, 0), (1, 8), (3, 6)\},~~
\{(0, 10), (1, 0), (3, 8)\},~~$

$\{(0, 0), (0, 3), (1, 10), (2, 6), (3, 7)\},~
\{(0, 10), (1, 7), (2, 5), (2, 6), (3, 0)\},~~
\{(0, 0), (1, 3), (3, 10)\},$

$\{(0, 1), (0, 9), (1, 1), (3, 0), (3, 6)\},\hspace{0.25cm}
\{(0, 7), (1, 0), (1, 9), (3, 0), (3, 10)\},~~
\{(0, 0), (0, 5), (1, 9), (3, 1)\}.$
\end{quote}}

\noindent13.~$(n,m)=(7,27):$~
{\scriptsize$$\begin{array}{lllll}
\hspace{-0.3cm}\{(0, 4), (1, 0), (3, 9)\},~~~~
\{(0, 2), (0, 6), (3, 0)\},\hspace{0.68cm}
\{(0, 0), (0, 11), (2, 2)\},&
\{(0, 0), (0, 5), (2, 11)\},\\
\hspace{-0.3cm}\{(0, 0), (0, 8), (3, 0)\},~~~~
\{(0, 1), (1, 0), (3, 8)\},\hspace{0.68cm}
\{(0, 0), (1, 3), (3, 13)\},&
\{(0, 6), (1, 0), (3, 12)\},\\
\hspace{-0.3cm}\{(0, 0), (0, 7), (3, 4)\},~~~~
\{(0, 0), (1, 10), (3, 3)\},\hspace{0.55cm}
\{(0, 0), (1, 8), (3, 11)\},&
\{(0, 11), (1, 0), (3, 13)\},\\
\hspace{-0.3cm}\{(0, 0), (1, 9), (3, 9)\},~~~~
\{(0, 0), (0, 13), (1, 11)\},\hspace{0.4cm}
\{(0, 4), (0, 10), (2, 0)\},&
\{(0, 0), (1, 5), (3, 12)\},\\
\hspace{-0.3cm}\{(0, 0), (1, 13), (2, 4), (3, 1), (3, 10)\},~~
\{(0, 10), (0, 13), (2, 8), (3, 0), (3, 1)\},&
\{(0, 0), (0, 12), (1, 7), (3, 8)\},\\
\hspace{-0.3cm}\{(0, 1), (0, 11), (1, 13), (2, 0), (3, 0)\},~~
\{(0, 10), (0, 12), (1, 0), (2, 4), (3, 7)\}.
\end{array}$$}

\noindent14.~$(n,m)=(9,9):$~
{\scriptsize \begin{quote}
$\{(0, 0), (1, 0), (4, 4)\},~~
\{(0, 4), (2, 0), (4, 0)\},~~
\{(0, 0), (3, 0), (4, 3)\},~~
\{(0, 1), (1, 0), (3, 4)\},~~$

$\{(0, 0), (4, 0), (4, 2)\},~~
\{(0, 0), (0, 3), (4, 1)\},~~
\{(0, 0), (1, 4), (2, 1)\},~~
\{(0, 2), (3, 0), (3, 4)\},~~$

$\{(0, 1), (2, 3), (4, 0)\},~~
\{(0, 4), (3, 0), (3, 1)\},~~
\{(0, 3), (1, 1), (2, 2), (3, 4), (4, 0)\}.$
\end{quote}}

\noindent15.~$(n,m)=(9,11):$~
{\scriptsize \begin{quote}
$\{(0, 0), (1, 0), (4, 5)\},~~
\{(0, 5), (2, 0), (4, 0)\},~~
\{(0, 0), (3, 0), (4, 4)\},~~
\{(0, 1), (1, 0), (3, 5)\},~~$

$\{(0, 0), (4, 0), (4, 3)\},~~
\{(0, 4), (1, 1), (4, 0)\},~~
\{(0, 5), (3, 0), (3, 1)\},~~
\{(0, 2), (3, 4), (4, 0)\},~~$

$\{(0, 2), (1, 5), (2, 0)\},~~
\{(0, 4), (2, 0), (4, 1)\},~~
\{(0, 2), (3, 0), (3, 5)\},~~
\{(0, 0), (0, 2), (2, 4), (4, 1)\},~~$

$\{(0, 3), (2, 2), (3, 0), (3, 4), (4, 5)\}.$
\end{quote}}\vspace{0.5cm}

\noindent16.~$(n,m)=(9,13):$~
{\scriptsize \begin{quote}
$\{(0, 5), (2, 0), (4, 0)\},~~
\{(0, 4), (3, 0), (4, 0)\},~~
\{(0, 1), (1, 6), (4, 0)\},~~
\{(0, 0), (1, 6), (3, 0)\},~~$

$\{(0, 3), (3, 0), (4, 1)\},~~
\{(0, 0), (2, 2), (3, 5)\},~~
\{(0, 5), (1, 0), (4, 6)\},~~
\{(0, 0), (1, 4), (3, 2)\},~~$

$\{(0, 3), (1, 0), (1, 5)\},~~
\{(0, 4), (2, 0), (2, 1)\},~~
\{(0, 0), (0, 4), (4, 6)\},~~
\{(0, 0), (2, 1), (4, 4)\},~~$

$\{(0, 1), (2, 0), (4, 4), (4, 6)\},~~
\{(0, 6), (1, 0), (3, 5), (4, 3)\},~~
\{(0, 0), (0, 6), (3, 1), (3, 4), (4, 0)\}.$
\end{quote}}

\noindent17.~$(n,m)=(9,15):$~
{\scriptsize \begin{quote}
$\{(0, 5), (1, 1), (4, 0)\},\hspace{0.34cm}
\{(0, 1), (0, 4), (4, 0)\},\hspace{0.34cm}
\{(0, 7), (2, 0), (3, 1)\},\hspace{0.3cm}
\{(0, 2), (1, 6), (3, 0)\},~~$

$\{(0, 1), (2, 0), (2, 7)\},\hspace{0.34cm}
\{(0, 0), (3, 2), (4, 7)\},\hspace{0.34cm}
\{(0, 0), (1, 2), (4, 5)\},\hspace{0.3cm}
\{(0, 0), (3, 6), (4, 1)\},~~$

$\{(0, 4), (1, 2), (3, 0), (4, 6), (4, 7)\},~~
\{(0, 0), (0, 2), (2, 5), (4, 0), (4, 6)\},~~~
\{(0, 0), (1, 3), (3, 7)\},$

$\{(0, 3), (0, 7), (2, 3), (3, 0), (4, 0)\},~~
\{(0, 7), (1, 0), (3, 2), (3, 7), (4, 1)\}.$
\end{quote}}

\noindent18.~$(n,m)=(9,17):$~
{\scriptsize \begin{quote}
$\{(0, 5), (2, 1), (4, 0)\},~~
\{(0, 0), (2, 8), (3, 2)\},~~
\{(0, 0), (2, 7), (4, 2)\},~~
\{(0, 0), (1, 8), (2, 3)\},$

$\{(0, 0), (3, 7), (4, 0)\},~~
\{(0, 0), (1, 7), (4, 3)\},~~
\{(0, 0), (3, 1), (3, 8)\},~~
\{(0, 0), (0, 3), (4, 7)\},$
\end{quote}
\begin{quote}
$\{(0, 0), (0, 8), (2, 2), (3, 0)\},~~
\{(0, 8), (1, 0), (2, 1), (3, 6), (4, 6)\},~~
\{(0, 3), (0, 8), (2, 0), (4, 0), (4, 4)\},$

$\{(0, 1), (1, 7), (4, 0), (4, 6)\},~~
\{(0, 0), (1, 3), (3, 4), (4, 6), (4, 8)\},~~
\{(0, 6), (0, 7), (1, 3), (3, 1), (4, 0)\}.$
\end{quote}}

\noindent19.~$(n,m)=(9,19):$~
{\scriptsize \begin{quote}
$\{(0, 9), (3, 0), (4, 0)\},~~
\{(0, 0), (2, 0), (4, 9)\},~~
\{(0, 8), (1, 0), (4, 0)\},~~
\{(0, 0), (4, 0), (4, 8)\},~~$

$\{(0, 1), (0, 7), (4, 0)\},~~
\{(0, 6), (1, 1), (4, 0)\},~~
\{(0, 9), (2, 0), (2, 1)\},~~
\{(0, 4), (2, 1), (4, 0)\},~~$

$\{(0, 7), (2, 1), (3, 0)\},~~
\{(0, 6), (3, 0), (4, 1)\},~~
\{(0, 4), (2, 0), (4, 1)\},~~
\{(0, 0), (0, 9), (3, 1)\},~~$

$\{(0, 2), (3, 0), (4, 8)\},~~
\{(0, 2), (2, 0), (4, 6)\},~~
\{(0, 0), (1, 7), (4, 2)\},~~
\{(0, 7), (1, 0), (3, 4)\},~~$

$\{(0, 0), (0, 7), (2, 2)\},~~
\{(0, 9), (1, 0), (3, 5)\},~~
\{(0, 7), (2, 0), (4, 8)\},~~
\{(0, 0), (3, 8), (4, 5)\},~~$

$\{(0, 0), (1, 4), (4, 7)\},~~
\{(0, 0), (1, 3), (3, 6)\},~~
\{(0, 2), (1, 0), (3, 7)\},~~
\{(0, 0), (0, 3), (1, 5), (1, 9)\},$

$\{(0, 6), (1, 0), (1, 2), (4, 4), (4, 9)\}.$
\end{quote}}

\noindent20.~$(n,m)=(9,21):$~
{\scriptsize \begin{quote}
$\{(0, 9), (1, 0), (4, 0)\},~~
\{(0, 0), (0, 8), (4, 0)\},~~
\{(0, 10), (3, 0), (4, 0)\},~~
\{(0, 0), (2, 0), (4, 10)\},~~$

$\{(0, 7), (1, 1), (4, 0)\},~~
\{(0, 6), (2, 1), (4, 0)\},~~
\{(0, 10), (2, 0), (2, 1)\},~~
\{(0, 1), (4, 0), (4, 10)\},~~$

$\{(0, 8), (2, 1), (3, 0)\},~~
\{(0, 5), (3, 0), (4, 1)\},~~
\{(0, 10), (1, 0), (3, 1)\},\hspace{0.2cm}
\{(0, 0), (0, 6), (4, 1)\},~~$

$\{(0, 4), (1, 0), (4, 1)\},~~
\{(0, 2), (2, 8), (4, 0)\},~~
\{(0, 2), (0, 7), (3, 0)\},\hspace{0.33cm}
\{(0, 8), (1, 0), (3, 2)\},~~$

$\{(0, 3), (3, 0), (4, 6)\},~~
\{(0, 0), (3, 4), (4, 8)\},~~
\{(0, 2), (1, 0), (1, 9)\},\hspace{0.33cm}
\{(0, 4), (2, 0), (4, 8)\},~~$

$\{(0, 0), (1, 8), (3, 2)\},~~
\{(0, 3), (2, 0), (3, 9)\},~~
\{(0, 0), (2, 5), (2, 9)\},\hspace{0.33cm}
\{(0, 0), (2, 4), (4, 7)\},~~$

$\{(0, 2), (1, 7), (2, 0)\},~~
\{(0, 0), (1, 2), (3, 9)\},~~
\{(0, 4), (3, 0), (3, 7), (4, 10)\},~~$

$\{(0, 5), (1, 0), (1, 2), (4, 7), (4, 10)\}.$
\end{quote}}

\noindent21.~$(n,m)=(9,23):$~
{\scriptsize$$\begin{array}{lllll}
\{(0, 0), (1, 2), (4, 7)\},&
\{(0, 5), (2, 0), (4, 1)\},&
\{(0, 0), (2, 0), (4, 11)\},&
\{(0, 0), (4, 0), (4, 10)\},\\
\{(0, 1), (0, 9), (4, 0)\},&
\{(0, 8), (1, 1), (4, 0)\},&
\{(0, 7), (2, 1), (4, 0)\},&
\{(0, 10), (2, 1), (3, 0)\},\\
\{(0, 6), (4, 0), (4, 1)\},&
\{(0, 0), (1, 4), (3, 8)\},&
\{(0, 10), (2, 0), (3, 1)\},&
\{(0, 10), (1, 0), (4, 0)\},\\
\{(0, 0), (1, 8), (4, 1)\},&
\{(0, 4), (1, 0), (4, 1)\},&
\{(0, 2), (0, 8), (3, 0)\},&
\{(0, 2), (1, 11), (2, 0)\},\\
\{(0, 6), (3, 0), (3, 2)\},&
\{(0, 0), (2, 7), (4, 3)\},&
\{(0, 3), (1, 0), (4, 7)\},&
\{(0, 11), (3, 0), (4, 0)\},\\
\{(0, 0), (2, 2), (3, 9)\},&
\{(0, 0), (3, 2), (4, 8)\},&
\{(0, 0), (0, 4), (3, 10)\},&
\{(0, 0), (1, 10), (3, 3)\},\\
\{(0, 0), (0, 9), (1, 3)\},&
\{(0, 2), (2, 11), (4, 0)\},&
\{(0, 0), (0, 11), (2, 3)\},&
\{(0, 0), (0, 5), (2, 10)\},\\
\{(0, 3), (2, 0), (4, 8)\},&
\{(0, 5), (3, 0), (4, 11)\},&
\{(0, 8), (1, 0), (1, 3)\},&
\{(0, 2), (0, 9), (1, 0), (3, 6), (4, 11)\}.
\end{array}$$}

\noindent22.~$(n,m)=(9,29):$~
{\scriptsize$$\begin{array}{lllll}
\hspace{0.7cm}\{(0, 6), (1, 6), (4, 0)\},&
\{(0, 0), (1, 1), (3, 7)\},&
\{(0, 7), (2, 0), (4, 8)\},\hspace{0.53cm}
\{(0, 0), (1, 8), (3, 6)\},\\
\hspace{0.7cm}\{(0, 0), (2, 0), (4, 7)\},&
\{(0, 0), (1, 9), (3, 0)\},&
\{(0, 5), (1, 2), (4, 0)\},\hspace{0.53cm}
\{(0, 0), (1, 6), (3, 1)\},\\
\hspace{0.7cm}\{(0, 9), (1, 3), (4, 0)\},&
\{(0, 0), (0, 6), (4, 11)\},&
\{(0, 7), (0, 14), (3, 0)\},\hspace{0.4cm}
\{(0, 3), (0, 11), (4, 0)\},\\
\hspace{0.7cm}\{(0, 0), (0, 9), (3, 13)\},&
\{(0, 0), (0, 10), (4, 9)\},&
\{(0, 0), (0, 12), (4, 8)\},\hspace{0.4cm}
\{(0, 0), (0, 11), (2, 14)\},\\
\hspace{0,7cm}\{(0, 0), (0, 13), (3, 12)\},&
\{(0, 0), (0, 14), (4, 4)\},&
\{(0, 14), (1, 0), (3, 1)\},\hspace{0.4cm}
\{(0, 14), (1, 1), (2, 0)\},\\
\hspace{0.7cm}\{(0, 12), (1, 0), (2, 14)\},&
\{(0, 7), (1, 0), (3, 12)\},&
\{(0, 0), (1, 10), (4, 2)\},\hspace{0.4cm}
\{(0, 10), (1, 0), (3, 13)\},\\
\hspace{0.7cm}\{(0, 12), (1, 3), (3, 0)\},&
\{(0, 0), (0, 5), (3, 14)\},&
\{(0, 11), (1, 6), (3, 0)\},\hspace{0.4cm}
\{(0, 10), (1, 13), (3, 0)\},\\
\hspace{0.7cm}\{(0, 13), (0, 14), (4, 0)\},&
\{(0, 0), (0, 2), (4, 14)\},&
\{(0, 8), (0, 12), (4, 0)\},\hspace{0.4cm}
\{(0, 0), (1, 11), (3, 10)\},\\
\hspace{0.7cm}\{(0, 8), (2, 0), (4, 11)\},&
\{(0, 4), (1, 0), (2, 13)\},&
\{(0, 12), (2, 0), (4, 5)\},\hspace{0.4cm}
\{(0, 9), (1, 11), (3, 0),(4, 7)\},\\
\hspace{0.7cm}\{(0, 0), (0, 3), (4, 13)\},&
\{(0, 5), (1, 10), (3, 0)\},&
\{(0, 11), (1, 0), (4, 11)\},\hspace{0.25cm}
\{(0, 2), (1, 0), (2, 12), (3, 14), (4, 8)\}.
\end{array}$$}

\noindent23.~$(n,m)=(11,11)$,
{\scriptsize$$\begin{array}{lllll}
\hspace{-1.3cm}\{(0,0),(1,0),(5,5)\},&
\{(0,5),(3,0),(5,0)\},&
\{(0,0),(3,0),(5,4)\},&
\{(0,4),(1,0),(5,0)\},\\
\hspace{-1.3cm}\{(0,0),(1,5),(5,0)\},&
\{(0,1),(0,3),(5,0)\},&
\{(0,1),(0,4),(4,0)\},&
\{(0,1),(2,0),(5,4)\},\\
\hspace{-1.3cm}\{(0,1),(1,0),(4,5)\},&
\{(0,1),(3,0),(5,3)\},&
\{(0,0),(1,3),(5,1)\},&
\{(0,5),(1,2),(2,0)\},\\
\hspace{-1.3cm}\{(0,2),(2,0),(4,5)\},&
\{(0,0),(1,1),(4,2)\},&
\{(0,2),(2,4),(5,0)\},&
\{(0,2),(3,0),(3,4),(3,5)\}\\
\hspace{-1.3cm}\{(0,4),(2,0),(4,1)\},&
\{(0,3),(3,0),(4,4)\},&
\{(0,3),(1,5),(2,0)\}.
\end{array}$$}

\noindent24.~$(n,m)=(11,13):$~
{\scriptsize \begin{quote}
$\{(0, 6), (4, 0), (5, 0)\},~~
\{(0, 0), (2, 0), (5, 6)\},~~
\{(0, 5), (2, 0), (5, 0)\},~~
\{(0, 0), (4, 0), (5, 5)\},~~$

$\{(0, 0), (4, 6), (5, 0)\},~~
\{(0, 1), (5, 0), (5, 5)\},~~
\{(0, 4), (1, 1), (5, 0)\},~~
\{(0, 6), (1, 1), (3, 0)\},~~$

$\{(0, 4), (3, 1), (4, 0)\},~~
\{(0, 5), (1, 1), (4, 0)\},~~
\{(0, 3), (3, 6), (5, 0)\},~~
\{(0, 2), (4, 0), (5, 4)\},~~$

$\{(0, 4), (3, 0), (4, 1)\},~~
\{(0, 0), (1, 2), (1, 6)\},~~
\{(0, 0), (2, 2), (3, 5)\},~~
\{(0, 0), (2, 1), (4, 5), (5, 3)\},~~$

$\{(0, 0), (0, 1), (2, 6), (4, 2), (4, 4)\},~~
\{(0, 3), (2, 0), (2, 6), (5, 1), (5, 4)\}.$
\end{quote}}

\noindent25.~$(n,m)=(11,15):$~
{\scriptsize \begin{quote}
$\{(0, 5), (0, 7), (4, 0)\},~~
\{(0, 4), (0, 5), (5, 0)\},~~
\{(0, 3), (0, 6), (3, 0)\},~~~
\{(0, 0), (0, 7), (5, 4)\},~~$

$\{(0, 7), (1, 0), (3, 5)\},~~
\{(0, 6), (1, 0), (4, 2)\},~~
\{(0, 3), (1, 0), (5, 3)\},~~~
\{(0, 0), (1, 7), (5, 5)\},~~$

$\{(0, 0), (1, 6), (4, 1)\},~~
\{(0, 4), (1, 0), (5, 7)\},~~
\{(0, 2), (1, 0), (4, 6)\},~~~
\{(0, 0), (2, 3), (5, 6)\},~~$

$\{(0, 1), (0, 7), (2, 7), (4, 1), (5, 0)\},~
\{(0, 3), (0, 7), (2, 0), (2, 5), (5, 1) \},~
\{(0, 0), (1, 4), (3, 5), (5, 1)\},$

$\{(0, 5), (1, 5), (2, 0), (4, 4), (5, 7)\},~
\{(0, 0), (1, 1), (3, 0), (4, 5), (5, 7)\}.$
\end{quote}}

\noindent26.~$(n,m)=(11,17)$,
{\scriptsize\begin{quote}
$\{(0, 0), (1, 0), (4, 6)\},~~
\{(0, 0), (1, 2), (1, 8)\},~~
\{(0, 3), (4, 8), (5, 0)\},~~
\{(0, 0), (2, 8), (3, 1)\},~~$

$\{(0, 3), (3, 0), (3, 8)\},~~
\{(0, 0), (1, 7), (5, 1)\},~~
\{(0, 1), (2, 6), (4, 0)\},~~
\{(0, 8), (2, 0), (5, 3)\},~~$

$\{(0, 4), (2, 8), (5, 0)\},~~
\{(0, 7), (4, 3), (5, 0)\},~~
\{(0, 3), (1, 6), (4, 0)\},~~
\{(0, 7), (1, 1), (4, 0)\},~~$

$\{(0, 7), (2, 0), (2, 2)\},~~
\{(0, 8), (2, 5), (5, 0)\},~~
\{(0, 6), (2, 2), (5, 0)\},~~
\{(0, 2), (1, 0), (4, 0), (4, 4)\},~~$

$\{(0, 0), (0, 4), (2, 3), (5, 7)\},~~
\{(0, 4), (1, 0), (4, 7), (5, 2)\},~~
\{(0, 2), (2, 0), (4, 2), (5, 8)\},~~$

$\{(0, 0), (2, 1), (4, 1), (4, 8)\},~~
\{(0, 0), (0, 5), (5, 4), (5, 5)\},~~
\{(0, 5), (0, 8), (3, 1), (4, 0)\},~~$

$\{(0, 0), (3, 2), (4, 7), (5, 8)\}.$
\end{quote}}

\noindent27.~$(n,m)=(11,21):$~
{\scriptsize$$\begin{array}{lllll}
\hspace{-1.7cm}\{(0, 2), (0, 8), (5, 0)\},\hspace{0.3cm}
\{(0, 0), (0, 7), (5, 7)\},\hspace{0.3cm}
\{(0, 0), (0, 8), (4, 9)\},\hspace{0.4cm}
\{(0, 5), (0, 10), (4, 0)\},\\
\hspace{-1.7cm}\{(0, 3), (1, 0), (5, 7)\},\hspace{0.3cm}
\{(0, 2), (1, 0), (2, 9)\},\hspace{0.3cm}
\{(0, 1), (1, 0), (3, 8)\},\hspace{0.4cm}
\{(0, 4), (1, 0), (4, 8)\},\\
\hspace{-1.7cm}\{(0, 0), (1, 0), (4, 5)\},\hspace{0.3cm}
\{(0, 0), (1, 6), (3, 3)\},\hspace{0.3cm}
\{(0, 6), (2, 1), (5, 0)\},\hspace{0.4cm}
\{(0, 8), (2, 2), (4, 0)\},\\
\hspace{-1.7cm}\{(0, 3), (2, 2), (5, 0)\},\hspace{0.3cm}
\{(0, 0), (2, 9), (5, 5)\},\hspace{0.3cm}
\{(0, 0), (1, 7), (4, 2)\},\hspace{0.4cm}
\{(0, 0), (1, 8), (4, 0)\},\\
\hspace{-1.7cm}\{(0, 7), (2, 0), (4, 6)\},\hspace{0.3cm}
\{(0, 2), (1, 6), (4, 0)\},\hspace{0.3cm}
\{(0, 1), (2, 3), (5, 0)\},\hspace{0.4cm}
\{(0, 8), (2, 0), (4, 4)\},\\
\hspace{-1.7cm}\{(0, 7), (1, 0), (5, 8)\},\hspace{0.3cm}
\{(0, 6), (1, 0), (4, 9)\},\hspace{0.3cm}
\{(0, 8), (1, 0), (3, 10)\},\hspace{0.26cm}
\{(0, 5), (1, 0), (2, 10)\},\\
\hspace{-1.7cm}\{(0, 0), (0, 3), (5, 9)\},\hspace{0.3cm}
\{(0, 9), (0, 10), (5, 0)\},\hspace{0.17cm}
\{(0, 0), (0, 2), (5, 10)\},\hspace{0.26cm}
\{(0, 0), (0, 4), (4, 10)\},\\
\hspace{-1.7cm}\{(0, 9), (1, 0), (4, 0)\},\hspace{0.3cm}
\{(0, 9), (2, 0), (5, 4)\},\hspace{0.3cm}
\{(0, 0), (0, 9), (3, 10)\},\hspace{0.28cm}
\{(0, 0), (0, 10), (5, 3)\},\\
\hspace{-1.7cm}\{(0, 10), (1, 0), (3, 0)\},\hspace{0.15cm}
\{(0, 0), (1, 5), (3, 6), (5, 2)\},\hspace{0.3cm}
\{(0, 7), (1, 9), (2, 10), (4, 0), (5, 3)\}.
\end{array}$$}

\noindent28.~$(n,m)=(13,15):$~
{\scriptsize \begin{quote}
$\{(0, 7), (5, 0), (6, 0)\},~~
\{(0, 0), (2, 0), (6, 7)\},~~
\{(0, 6), (3, 0), (6, 0)\},~~
\{(0, 0), (4, 0), (6, 6)\},~~$

$\{(0, 0), (5, 0), (5, 7)\},~~
\{(0, 0), (6, 0), (6, 5)\},~~
\{(0, 1), (5, 7), (6, 0)\},~~
\{(0, 5), (1, 1), (6, 0)\},~~$

$\{(0, 1), (0, 7), (4, 0)\},~~
\{(0, 4), (3, 1),(6, 0)\},~~
\{(0, 6), (3, 1), (5, 0)\},~~
\{(0, 7), (2, 1), (3, 0)\},~~$

$\{(0, 3), (6, 0), (6, 1)\},~~
\{(0, 6), (4, 0), (5, 1)\},~~
\{(0, 4), (3, 0), (5, 1)\},~~
\{(0, 0), (0, 2), (6, 4)\},~~$

$\{(0, 5), (2, 0), (4, 2)\},~~
\{(0, 5), (1, 0), (5, 1)\},~~
\{(0, 0), (0, 3), (3, 6)\},~~
\{(0, 2), (2, 0), (4, 5)\},~~$

$\{(0, 0), (1, 7), (5, 2)\},~~
\{(0, 0), (3, 1), (5, 5)\},~~
\{(0, 0), (1, 3), (2, 7)\},~~
\{(0, 0), (1, 6), (3, 2)\},~~$

$\{(0, 2), (1, 0), (5, 5)\},~~
\{(0, 2), (3, 6), (4, 0)\},~~
\{(0, 0), (4, 2), (4, 6)\},~~
\{(0, 3), (2, 0), (5, 7), (6, 4)\},$

$\{(0, 2), (1, 4), (3, 7), (5, 0), (6, 5)\}.$
\end{quote}}

\noindent29.~$(n,m)=(15,15):$~
{\scriptsize \begin{quote}
$\{(0, 7), (6, 0), (7, 0)\},~~
\{(0, 0), (2, 0), (7, 7)\},~~
\{(0, 6), (4, 0), (7, 0)\},~~
\{(0, 0), (4, 0), (7, 6)\},~~$

$\{(0, 0), (5, 0), (6, 7)\},~~
\{(0, 0), (6, 0), (7, 5)\},~~
\{(0, 0), (6, 6), (7, 0)\},~~
\{(0, 1), (2, 7), (7, 0)\},~~$

$\{(0, 1), (5, 0), (7, 5)\},~~
\{(0, 1), (1, 0), (5, 7)\},~~
\{(0, 1), (6, 0), (7, 4)\},~~
\{(0, 1), (2, 0), (6, 6)\},~~$

$\{(0, 1), (3, 0), (6, 5)\},~~
\{(0, 1), (4, 0), (7, 2)\},~~
\{(0, 0), (0, 1), (6, 3)\},~~
\{(0, 4), (6, 0), (7, 1)\},~~$

$\{(0, 6), (3, 0), (6, 1)\},~~
\{(0, 4), (4, 0), (6, 1)\},~~
\{(0, 6), (2, 0), (7, 1)\},~~
\{(0, 2), (7, 0), (7, 4)\},~~$

$\{(0, 6), (1, 2), (5, 0)\},~~
\{(0, 2), (5, 0), (5, 6)\},~~
\{(0, 6), (3, 2), (6, 0)\},~~
\{(0, 2), (2, 0), (2, 5)\},~~$

$\{(0, 0), (4, 5), (5, 2)\},~~
\{(0, 0), (3, 7), (5, 3)\},~~
\{(0, 0), (0, 3), (2, 7)\},~~
\{(0, 0), (1, 6), (6, 1)\},~~$

$\{(0, 0), (4, 2), (5, 5)\},~~
\{(0, 3), (4, 7), (5, 0)\},~~
\{(0, 3), (3, 0), (3, 7)\},~~
\{(0, 5), (2, 0), (4, 2)\},~~$

$\{(0, 7), (4, 0), (4, 2)\},~~
\{(0, 7), (3, 0), (7, 3)\},~~
\{(0, 5), (1, 7), (2, 2), (3, 0), (6, 3)\}.$
\end{quote}}

\section{Generalized $([n]\times M,\{3,4,5\},1)$-PDFs}\label{app-B}
The base blocks of a generalized $([n]\times M,\{3,4,5\},1)$-PDF are listed below:\\

\noindent1.~$n=15, M=\{0,\pm1,\pm7,\pm8\}:$
{\scriptsize \begin{quote}
$\{(0, 0), (0, 1), (7, 1)\},~~
\{(0, 0), (1, 0), (7, 8)\},~~
\{(0, 0), (0, 8), (6, 1)\},~~
\{(0, 0), (0, 7), (5, 8)\},~~$

$\{(0, 0), (1, 8), (6, 7)\},~~
\{(0, 0), (1, 7), (7, 7)\},~~
\{(0, 0), (1, 1), (4, 8)\},~~
\{(0, 1), (1, 0), (4, 8)\},~~$

$\{(0, 7), (2, 0), (4, 8)\},~~
\{(0, 1), (5, 8), (6, 0)\},~~
\{(0, 8), (4, 0), (7, 1)\},~~
\{(0, 1), (3, 0), (7, 0)\},~~$

$\{(0, 8), (3, 0), (5, 1)\},~~
\{(0, 8), (2, 7), (6, 0)\},~~
\{(0, 8), (2, 0), (4, 7), (5, 0), (7, 0)\}$.
\end{quote}}

\noindent2.~$n=17, M=\{0,\pm1,\pm7,\pm8\}:$
{\scriptsize \begin{quote}
$\{(0, 0), (0, 1), (8, 1)\},~~
\{(0, 0), (1, 0), (8, 8)\},~~
\{(0, 0), (0, 8), (7, 1)\},~~
\{(0, 0), (1, 8), (8, 7)\},~~$

$\{(0, 0), (0, 7), (6, 8)\},~~
\{(0, 0), (2, 7), (5, 8)\},~~
\{(0, 1), (3, 0), (6, 8)\},~~
\{(0, 1), (1, 0), (5, 8)\},~~$

$\{(0, 1), (3, 8), (8, 0)\},~~
\{(0, 8), (2, 1), (7, 0)\},~~
\{(0, 8), (4, 7), (8, 0)\},~~
\{(0,~7), (1,~0), (6,~0)\},~~$

$\{(0,~7), (1,~8), (3,~0), (5,~8), (7,~7)\},~~
\{(0,~8), (3,~8), (4,~0), (6,~0), (8,~1)\}$.
\end{quote}}

\noindent3.~$n=21, M=\{0,\pm1,\pm7,\pm8\}:$
{\scriptsize \begin{quote}
$\{(0, 0), (0, 1), (10, 1)\},\hspace{0.2cm}
\{(0, 0), (0, 8), (9, 1)\},\hspace{0.2cm}
\{(0, 0), (1, 0), (10, 8)\},\hspace{0.15cm}
\{(0, 0), (1, 8), (10, 7)\},$

$\{(0, 0), (0, 7), (8, 8)\},\hspace{0.35cm}
\{(0, 0), (1, 7), (9, 7)\},\hspace{0.2cm}
\{(0, 0), (1, 1), (7, 8)\},\hspace{0.3cm}
\{(0, 0), (2, 0), (5, 8)\},$

$\{(0, 0), (2, 8), (7, 7)\},\hspace{0.35cm}
\{(0, 0), (2, 1), (6, 8)\},\hspace{0.2cm}
\{(0, 0), (3, 1), (7, 1)\},\hspace{0.3cm}
\{(0, 0), (4, 8), (6, 1)\},$

$\{(0, 8), (2, 0), (10, 7)\},\hspace{0.2cm}
\{(0, 1), (4, 0), (9, 1)\},\hspace{0.2cm}
\{(0, 8), (5, 0), (7, 7)\},\hspace{0.3cm}
\{(0, 8), (3, 0), (8, 7)\},~~$

$\{(0, 8), (3, 7), (10, 0)\},\hspace{0.2cm}
\{(0, 8), (2, 7), (8, 0)\},\hspace{0.2cm}
\{(0, 8), (6, 7), (7, 0)\},\hspace{0.3cm}
\{(0, 7), (4, 0), (7, 7), (10, 0)\},~~$

$\{(0, 8), (1, 7), (5, 8), (6, 0), (9, 0)\}$.
\end{quote}}

\noindent4.~$n=27, M=\{0,\pm1,\pm7,\pm8\}:$
{\scriptsize$$\begin{array}{llllll}
\hspace{0.4cm}\{(0, 0), (0, 1), (13, 1)\},&
\{(0, 0), (1, 0), (13, 8)\},&
\{(0, 0), (0, 8), (12, 1)\},\hspace{0.4cm}
\{(0, 0), (0, 7), (11, 8)\},\\
\hspace{0.4cm}\{(0, 0), (1, 8), (12, 7)\},&
\{(0, 0), (1, 7), (13, 7)\},&
\{(0, 0), (1, 1), (10, 8)\},\hspace{0.4cm}
\{(0, 0), (2, 8), (11, 7)\},\\
\hspace{0.4cm}\{(0, 0), (2, 0), (9, 8)\},&
\{(0, 0), (2, 1), (10, 1)\},&
\{(0, 0), (2, 7), (8, 8)\},\hspace{0.55cm}
\{(0, 0), (5, 8), (6, 0)\},\\
\hspace{0.4cm}\{(0, 0), (3, 8), (9, 1)\},&
\{(0, 0), (3, 1), (6, 8)\},&
\{(0, 0), (3, 0), (10, 7)\},\hspace{0.4cm}
\{(0, 0), (4, 8), (7, 1)\},\\
\hspace{0.4cm}\{(0, 8), (5, 0), (13, 1)\},&
\{(0, 8), (1, 7), (10, 0)\},&
\{(0, 8), (4, 1), (8, 0)\},\hspace{0.55cm}
\{(0, 8), (4, 8), (11, 0)\},\\
\hspace{0.4cm}\{(0, 8), (5, 1), (12, 0)\},&
\{(0, 8), (3, 7), (13, 0)\},&
\{(0, 1), (6, 8), (12, 0)\},\hspace{0.4cm}
\{(0, 7), (2, 0), (10, 7), (11, 0)\},\\
\hspace{0.4cm}\{(0, 1), (5, 0), (10, 0)\},&
\{(0, 8), (4, 0), (8, 1)\},&
\{(0, 1), (5, 8), (8, 0)\},
\hspace{0.55cm}\{(0, 1), (4, 8), (6, 0), (11, 1), (13, 0)\}.
\end{array}$$}

\section{$(n,\{3,4\},1)$-PDPs with leave $[h]^r$}\label{app-C}
The base blocks of an $(n,\{3,4\},1)$-PDP with leave $[h]^r$ are listed below:\\

\noindent~1.
$(n,h,r)=(25,7,3)$:\vspace{0.22cm}
{\scriptsize$$\begin{array}{llllllll}
\hspace{-7.25cm}\{0, 2, 7\}, &\{0, 8, 12\}, &\{0, 10, 11\}.
\end{array}$$}
\noindent2.
$(n,h,r)=(31,7,3)$:\vspace{0.22cm}
{\scriptsize$$\begin{array}{llllllll}
\hspace{-7.2cm}\{0, 1, 15\}, &\{0, 4, 11\}, &\{0, 5, 13\}, &\{0, 10, 12\}.
\end{array}$$}
\noindent3.
$(n,h,r)=(33,9,3)$:\vspace{0.22cm}
{\scriptsize$$\begin{array}{llllllll}
\hspace{-7.2cm}\{0, 1, 8\}, &\{0, 5, 16\}, &\{0, 10, 14\}, &\{0, 13, 15\}.
\end{array}$$}
\noindent4.
$(n,h,r)=(37,7,3)$:\vspace{0.22cm}
{\scriptsize$$\begin{array}{llllllll}
\hspace{-7.2cm}\{0, 1, 14\}, &\{0, 10, 18\}, &\{0, 15, 17\}, &\{ 0, 4, 11, 16\}.
\end{array}$$}
\noindent5.
$(n,h,r)=(39,9,3)$:\vspace{0.22cm}
{\scriptsize$$\begin{array}{llllllll}
\hspace{-6.3cm}\{0, 1, 11\}, &\{0, 2, 18\}, &\{0, 5, 19\}, &\{0, 7, 15\}, &\{0, 13, 17\}.
\end{array}$$}
\noindent6.
$(n,h,r)=(41,5,3)$:\vspace{0.22cm}
{\scriptsize$$\begin{array}{llllllll}
\hspace{-5.8cm}\{0, 5, 19\}, &\{0, 8, 18\}, &\{0, 11, 12\}, &\{0, 15, 17\}, &\{0, 4, 13, 20\}.
\end{array}$$}
\noindent7.
$(n,h,r)=(43,7,3)$:\vspace{0.22cm}
{\scriptsize$$\begin{array}{llllllll}
\hspace{-5.7cm}\{0, 2, 13\}, &\{0, 12, 20\}, &\{ 0, 16, 21\}, &\{0, 10, 17\}, &\{0, 1, 15, 19\}.
\end{array}$$}
\noindent8.
$(n,h,r)=(45,9,3)$:\vspace{0.22cm}
{\scriptsize$$\begin{array}{llllllll}
\hspace{-5.7cm}\{0, 4, 20\}, &\{0, 5, 18\}, &\{0, 11, 21\}, &\{0, 17, 19\}, &\{0, 14, 15, 22\}.
\end{array}$$}
\noindent9.
$(n,h,r)=(47,5,3)$:\vspace{0.22cm}
{\scriptsize$$\begin{array}{llllllll}
\hspace{-4.4cm}\{0, 1, 21\}, &\{0, 4, 18\}, &\{0, 5, 16\}, &\{0, 8, 23\}, &\{0, 9, 22\}, &\{0, 2, 12, 19\}.
\end{array}$$}
\noindent10.
$(n,h,r)=(51,9,3)$:\vspace{0.22cm}
{\scriptsize$$\begin{array}{llllllll}
\hspace{-4.2cm}\{0, 4, 17\}, &\{0, 5, 19\}, &\{0, 7, 18\}, &\{0, 8, 24\}, &\{0, 10, 25\}, &\{0, 2, 22, 23\}.
\end{array}$$}
\noindent11.
$(n,h,r)=(51,15,3)$:\vspace{0.22cm}
{\scriptsize$$\begin{array}{llllllll}
\hspace{-5.6cm}\{0, 1, 8\}, &\{0, 10, 23\}, &\{0, 11, 25\}, &\{0, 16, 20\}, &\{0, 2, 19, 24\}.
\end{array}$$}
\noindent12.
$(n,h,r)=(57,15,3)$:
{\scriptsize$$\begin{array}{llllllll}
\hspace{-2.5cm}\{0, 5, 28\}, &\{0, 8, 19\}, &\{0, 22, 24\}, &\{0, 25, 26\}, &\{0, 16, 20\}, &\{0, 14, 27\}, &\{0, 7, 17\}.
\end{array}$$}
\noindent13.
$(n,h,r)=(57,9,3)$:
{\scriptsize$$\begin{array}{llllllll}
\hspace{-1cm}\{0, 2, 26\}, &\{0, 4, 20\}, &\{0, 5, 28\}, &\{0, 8, 25\}, &\{0, 11, 21\}, &\{0, 13, 27\}, &\{0, 15, 22\},  &\{0, 18, 19\}.
\end{array}$$}
\noindent14.
$(n,h,r)=(59,11,3)$:
{\scriptsize$$\begin{array}{llllllll}
\hspace{-1.0cm}\{0, 10, 29\}, &\{0, 11, 28\}, &\{0, 4, 27\}, &\{0, 21, 26\}, &\{0, 7, 25\}, &\{0, 20, 22\}, &\{0, 8, 24\}, &\{0, 1, 14\}.
\end{array}$$}
\noindent15.
$(n,h,r)=(61,7,3)$:
{\scriptsize$$\begin{array}{llllllll}
\hspace{-0.4cm}\{0, 19, 29\}, &\{0, 4, 24\}, &\{0, 1, 27\}, &\{0, 14, 22\}, &\{0, 18, 25\}, &\{0, 11, 23\}, &\{0, 16, 21\}, &\{0, 2, 17,  30\}.
\end{array}$$}
\noindent16.
$(n,h,r)=(61,19,3)$:
{\scriptsize$$\begin{array}{llllllll}
\hspace{-2.7cm}\{0, 16, 29\}, &\{0, 5, 19\}, &\{0, 8, 25\}, &\{0, 26, 28\}, &\{0, 1, 23\}, &\{0, 20, 30\}, &\{0, 7, 11\}.
\end{array}$$}
\noindent17.
$(n,h,r)=(63,15,3)$:
{\scriptsize$$\begin{array}{llllllll}
\hspace{-0.7cm}\{0, 7, 31\}, &\{0, 25, 29\}, &\{0, 26, 27\}, &\{0, 13, 23\}, &\{0, 8, 28\}, &\{0, 5, 22\}, &\{0, 19, 30\}, &\{0, 14, 16\}.
\end{array}$$}
\noindent18.
$(n,h,r)=(65,5,3)$:
{\scriptsize$$\begin{array}{llllllll}
\hspace{-0.0cm}\{0, 5, 27\}, &\{0, 1, 20\}, &\{0, 13, 25\}, &\{0, 14, 23\}, &\{0, 7, 24\}, &\{0, 15, 31\}, &\{0, 4, 30, 32\},&\{0, 11, 21, 29\}.
\end{array}$$}
\noindent19.
$(n,h,r)=(65,11,3)$:
{\scriptsize$$\begin{array}{lllllll}
\hspace{-1.4cm}\{0, 8, 22\}, &\{0, 10, 26\}, &\{0, 27, 29\}, &\{0, 13, 31\}, &\{0, 23, 30\}, &\{0, 4, 21, 32\}, &\{ 0, 5, 24, 25\}.
\end{array}$$}
\noindent20.
$(n,h,r)=(67,7,3)$:
{\scriptsize$$\begin{array}{lllllll}
\hspace{-1.1cm}\{0, 1, 27\}, &\{0, 23, 33\}, &\{0, 13, 31\}, &\{0, 8, 29\}, &\{0, 2, 17, 24\}, &\{0, 4, 20, 32\}, &\{0, 11, 25, 30\}.
\end{array}$$}
\noindent21.
$(n,h,r)=(67,19,3)$:
{\scriptsize$$\begin{array}{llllllll}
\hspace{-0.7cm}\{0, 22, 32\}, &\{0, 2, 31\}, &\{0, 8, 28\}, &\{0, 14, 33\}, &\{0, 25, 26\}, &\{0, 11, 16\}, &\{0, 7, 30\}, &\{0, 4, 17\}.
\end{array}$$}
\noindent22.
$(n,h,r)=(69,9,3)$:
{\scriptsize$$\begin{array}{llllllll}
\hspace{0.2cm}\{0, 4, 31\}, &\{0, 16, 29\}, &\{0, 22, 30\}, &\{0, 20, 25\}, &\{0, 18, 32\}, &\{0, 21, 28\}, &\{0, 11, 26\}, &\{0, 10, 33, 34\},\\
\hspace{0.2cm}\{0, 17, 19\}.&&&&&
\end{array}$$}
\noindent23.
$(n,h,r)=(69,21,3)$:
{\scriptsize$$\begin{array}{llllllll}
\hspace{-0.6cm}\{0, 26, 34\}, &\{0, 2, 31\}, &\{0, 11, 25\}, &\{0, 5, 33\}, &\{0, 10, 32\}, &\{0, 19, 23\}, &\{0, 13, 20\}, &\{0, 1, 17\}.
\end{array}$$}
\noindent24.
$(n,h,r)=(71,5,3)$:
{\scriptsize$$\begin{array}{llllllll}
\hspace{0.45cm}\{0, 9, 31\}, &\{0, 27, 28\}, &\{0, 15, 29\}, &\{0, 11, 30\}, &\{0, 10, 35\}, &\{0, 18, 23\}, &\{0, 17, 21, 33\}, &\{0, 2, 26, 34\}, \\
\hspace{0.45cm}\{0, 7, 20\}.&&&&&&
\end{array}$$}
\noindent25.
$(n,h,r)=(71,11,3)$:
{\scriptsize$$\begin{array}{llllllll}
\{0, 8, 33\}, &\{0, 5, 29\}, &\{0, 10, 26\}, &\{0, 22, 23\}, &\{0, 20,  34\}, &\{0, 27, 31\}, &\{0, 18, 35\}, &\{0, 2, 13, 32\},\\
\{0, 21, 28\}.&&&&&&
\end{array}$$}
\noindent26.
$(n,h,r)=(73,19,4)$:
{\scriptsize$$\begin{array}{lllllllll}
\hspace{-0.55cm}\{0, 2, 35\}, &\{0, 21, 34\}, &\{0, 26, 31\}, &\{0, 9, 27\}, &\{0, 7, 30\}, &\{0, 10, 29\}, &\{0, 3, 25\}, &\{0, 11, 17\},\\
\hspace{-0.55cm}\{0, 14, 15\}.
\end{array}$$}
\noindent27.
$(n,h,r)=(75,9,3)$:
{\scriptsize$$\begin{array}{llllllll}
\hspace{0.65cm}\{0, 8, 31\}, &\{0, 2, 28\}, &\{0, 22, 32\}, &\{0, 20, 27\}, &\{0, 18, 35\}, &\{0, 19, 33\},   &\{0, 21, 36, 37\},&\{0, 5, 30, 34\},\\
\hspace{0.65cm}\{0, 13, 24\}.
\end{array}$$}
\noindent28.
$(n,h,r)=(75,21,3)$:
{\scriptsize$$\begin{array}{llllllll}
\hspace{-0.15cm}\{0, 33, 37\}, &\{0, 25, 32\}, &\{0, 16, 35\}, &\{0, 29, 31\}, &\{0, 26, 36\}, &\{0, 14, 34\}, &\{0, 5, 13\}, &\{0, 11, 28\},\\
\hspace{-0.15cm}\{0, 22, 23\}.&&&&&&&
\end{array}$$}
\noindent29.
$(n,h,r)=(79,25,3)$:
{\scriptsize$$\begin{array}{lllllllll}
\hspace{-0.5cm}\{0, 28, 38\}, &\{0, 4, 35\}, &\{0, 7, 29\}, &\{0, 32, 37\}, &\{0, 23, 34\}, &\{0, 17, 25\}, &\{0, 26, 39\}, &\{0, 2, 16\},\\
\hspace{-0.5cm}\{0, 1, 20\}.
\end{array}$$}
\noindent30.
$(n,h,r)=(81,15,3)$:
{\scriptsize$$\begin{array}{llllllll}
\hspace{-0.35cm}\{0, 8, 40\}, &\{0, 14, 39\}, &\{0, 22, 38\}, &\{0, 24, 37\}, &\{0, 10, 30\}, &\{0, 17, 28\}, &\{0, 4, 27\}, &\{0, 7, 26\},\\
\hspace{-0.35cm}\{0, 5, 34\}, &\{0, 1, 36\}, &\{0, 31, 33\}.&&&&&
\end{array}$$}
\noindent31.
$(n,h,r)=(85,7,3)$:
{\scriptsize$$\begin{array}{llllllll}
\hspace{-0.8cm}\{0, 28, 40\}, &\{0, 20, 27\}, &\{0, 17, 38\}, &\{0, 8, 31\}, &\{0, 11, 15, 41\}, &\{0, 5, 24, 42\}, &\{0, 1, 34, 36\}, \\
\hspace{-0.8cm}\{0, 13, 29\}, &\{0, 22, 32\},  &\{0, 25, 39\}.
\end{array}$$}
\noindent32.
$(n,h,r)=(85,19,3)$:
{\scriptsize$$\begin{array}{llllllll}
\hspace{-0.25cm}\{0, 10, 41\}, &\{0, 35, 40\}, &\{0, 25, 42\}, &\{0, 38, 39\}, &\{0, 8, 37\}, &\{0, 28, 30\}, &\{0, 23, 34\}, &\{0, 13, 32\}, \\
\hspace{-0.25cm}\{0, 26, 33\}, &\{0, 14, 36\}, &\{0, 4, 20\}.&&&&&
\end{array}$$}
\noindent33.
$(n,h,r)=(87,15,3)$:
{\scriptsize$$\begin{array}{llllllll}
\{0,8,38\},&\{0, 17, 43\},&\{0, 37, 42\},&\{0, 13, 41\},&\{0, 24, 40\},&\{0, 20, 39\},&\{0, 14, 25\},&\{0, 22, 32\},\\
\{0, 23, 27\},&\{0, 7, 36\},&\{0, 2, 33\},&\{0, 34, 35\}.&&&&
\end{array}$$}
\noindent34.
$(n,h,r)=(89,5,3)$:
{\scriptsize$$\begin{array}{llllllll}
\hspace{-0.1cm}\{0, 4, 36\}, &\{0, 22, 41\}, &\{0, 21, 33\}, &\{0, 28, 30, 37\}, &\{0, 26, 34, 44\},   &\{0, 15, 38, 39\}, &\{0, 16, 29, 43\}, \\
\hspace{-0.1cm}\{0, 20, 31\},  &\{0, 5, 40\}, &\{0, 25, 42\}.&&&&
\end{array}$$}
\noindent35.
$(n,h,r)=(89,11,3)$:
{\scriptsize$$\begin{array}{llllllll}
\hspace{-0.95cm}\{0, 1, 33\}, &\{0, 4, 29\}, &\{0, 28, 42\},  &\{0, 17, 43\}, &\{0, 23, 41\}, &\{ 0, 31, 39, 44\}, &\{0, 10, 30, 37\},\\
\hspace{-0.95cm}\{0, 24, 35\}, &\{0, 19, 40\},  &\{0, 22, 38\}, &\{0, 34, 36\}.
\end{array}$$}
\noindent36.
$(n,h,r)=(91,7,3)$:
{\scriptsize$$\begin{array}{llllllll}
\hspace{-0.6cm}\{0, 1, 45\}, &\{0, 13, 33\},   &\{0, 25, 43\}, &\{0, 28, 35\},  &\{0, 11, 32, 40\}, &\{0, 12, 14, 38\}, &\{0, 15, 37, 42\},\\
\hspace{-0.6cm} \{0, 4, 34\}, &\{0, 19, 36\}, &\{0, 16, 39\},  &\{0, 31, 41\}.
\end{array}$$}
\noindent37.
$(n,h,r)=(91,19,3)$:
{\scriptsize$$\begin{array}{llllllll}
\hspace{-0.15cm}\{0, 28, 44\}, &\{0, 35, 45\}, &\{0, 20, 43\}, &\{0, 8, 42\}, &\{0, 36, 41\}, &\{0, 22, 39\}, &\{0, 7, 40\}, &\{0, 37, 38\},\\
\hspace{-0.15cm}\{0, 29, 31\}, &\{0, 19, 32\}, &\{0, 4, 30\}, &\{0, 11, 25\}.&&&&
\end{array}$$}
\noindent38.
$(n,h,r)=(93,9,3)$:
{\scriptsize$$\begin{array}{llllllll}
\hspace{0.35cm}\{0, 11, 40\}, &\{0, 4, 32\},&\{0, 22, 42\}, &\{0, 16, 30\}, &\{0, 7, 41\}, &\{0, 5, 36\},   &\{0, 1, 25, 46\},  &\{0, 33, 35, 43\}, \\
\hspace{0.35cm}\{0, 27, 44\}, &
\{0, 18, 37\}, &\{0, 26, 39\}, &\{0, 15, 38\}.
\end{array}$$}
\noindent39.
$(n,h,r)=(93,21,3)$:
{\scriptsize$$\begin{array}{llllllll}
\hspace{-0.35cm}\{0, 8, 46\}, &\{0, 4, 44\}, &\{0, 22, 35\},  &\{0, 23, 37\}, &\{0, 41, 43\}, &\{0, 1, 34\}, &\{0, 5, 36\}, &\{0, 11, 28\},  \\
\hspace{-0.35cm}\{0, 32, 39\}, &\{0, 20, 45\}, &\{0, 16, 42\}, &\{0, 19, 29\}.
\end{array}$$}
\noindent40.
$(n,h,r)=(95,5,3)$:
{\scriptsize$$\begin{array}{lllllllll}
\hspace{-0.6cm}\{0,7,42\}, &\{0, 15, 47\}, &\{0, 13, 31\}, &\{0, 23, 43\},  &\{ 0, 40, 41, 45\}, &\{0, 19, 27, 44\}, &\{0, 2, 28, 39\}, \\
\hspace{-0.6cm}\{0, 14, 36\}, &\{0, 30, 46\}, &\{0, 9, 38\}, &\{0, 21, 33\}, &\{0, 24, 34\}.
\end{array}$$}
\noindent41.
$(n,h,r)=(95,11,3)$:
{\scriptsize$$\begin{array}{lllllllll}
\hspace{-0.75cm} \{0, 2, 25\}, &\{0, 27, 37\}, &\{0, 28, 47\},  &\{0, 22, 42\}, &\{0, 5, 26, 44\}, &\{0, 1, 33, 41\},  &\{0, 17, 30, 46\},\\
\hspace{-0.75cm}\{ 0, 24, 35\},  &\{0, 34, 38\}, &\{0, 36, 43\}, &\{0, 14, 45\}.
\end{array}$$}
\noindent42.
$(n,h,r)=(97,25,3)$:
{\scriptsize$$\begin{array}{llllllll}
\hspace{-0.15cm}\{0, 34, 47\}, &\{0, 25, 44\}, &\{0, 37, 41\}, &\{0, 10, 45\}, &\{0, 29, 46\}, &\{0, 14, 22\}, &\{0, 5, 31\}, &\{0, 32, 43\},\\
\hspace{-0.15cm}\{0, 20, 48\}, &\{0, 40, 42\}, &\{0, 1, 39\}, &\{0, 7, 23\}.&&&&
\end{array}$$}
\noindent43.
$(n,h,r)=(99,9,3)$:
{\scriptsize$$\begin{array}{llllllll}
\hspace{-0.15cm}\{0, 14, 38\}, &\{0, 2, 37\}, &\{0, 8, 40\}, &\{0, 30, 49\}, &\{0, 23, 44, 48\}, &\{0, 15, 28, 46\}, &\{0, 7, 33, 43\}, \\
\hspace{-0.15cm}\{0, 41, 42\}, &\{0, 22, 39\},   &\{0, 27, 47\},  &\{ 0, 11, 16, 45\}.&&&
\end{array}$$}
\noindent44.
$(n,h,r)=(99,21,3)$:
{\scriptsize$$\begin{array}{llllllll}
\hspace{-0.25cm}\{0, 17, 49\}, &\{0, 39, 46\}, &\{0, 2, 40\}, &\{0, 4, 47\}, &\{0, 25, 48\}, &\{0, 37, 42\}, &\{0, 26, 45\}, &\{0, 13, 29\}, \\
\hspace{-0.25cm}\{0, 31, 41\}, &\{0, 33, 44\}, &\{0, 1, 35\}, &\{0, 8, 28\}, &\{0, 22, 36\}.&&&
\end{array}$$}
\noindent45.
$(n,h,r)=(103,25,3)$:
{\scriptsize$$\begin{array}{llllllll}
\{0, 45, 50\}, &\{0, 8, 46\}, &\{0, 1, 43\}, &\{0, 2, 25\}, &\{0, 20, 37\}, &\{0, 31, 47\}, &\{0, 4, 26\}, &\{0, 11, 51\}, \\
\{0, 29, 48\}, &\{0, 28, 41\}, &\{0, 35, 49\}, &\{0, 32, 39\}, &\{0, 34, 44\}.&&&
\end{array}$$}
\noindent46.
$(n,h,r)=(109,7,3)$:
{\scriptsize$$\begin{array}{lllllll}
\hspace{-0.1cm}\{0, 1, 40\}, &\{0, 27, 29\}, &\{0, 20, 51\}, &\{0, 11, 52\},  &\{0, 21, 28, 46\},&\{0, 10, 42, 54\}, &\{0, 15, 37, 53\}, \\
\hspace{-0.1cm} \{0, 24, 47\}, &\{0, 45, 49\}, &\{0, 26, 34\},  &\{0, 5, 35, 48\}, &\{0, 14, 33, 50\}.
\end{array}$$}
\noindent47.
$(n,h,r)=(109,19,3)$:
{\scriptsize$$\begin{array}{llllllll}
\hspace{-0.25cm}\{0, 20, 54\}, &\{0, 5, 53\}, &\{0, 2, 52\}, &\{0, 35, 49\}, &\{0, 13, 46\}, &\{0, 26, 45\}, &\{0, 23, 51\}, &\{0, 4, 36\},\\
\hspace{-0.25cm}\{0, 1, 44\}, &\{0, 25, 42\}, &\{0, 31, 41\}, &\{0, 16, 38\}, &\{0, 30, 37\}, &\{0, 11, 40\}, &\{0, 39, 47\}.&
\end{array}$$}
\noindent48.
$(n,h,r)=(113,5,3)$:
{\scriptsize$$\begin{array}{llllllll}
\hspace{-0.45cm}\{0, 1, 36\}, &\{0, 11, 38\}, &\{0, 24, 42\}, &\{0, 32, 49\}, &\{0, 4, 33, 52\}, &\{0, 30, 50, 55\}, &\{0, 16, 47, 56\}, \\
\hspace{-0.45cm}\{0, 23, 44\}, &\{0, 14, 51\}, &\{0, 39, 46\},  &\{0, 22, 34\},  &\{0, 13, 28, 54\}, &\{0, 8, 10, 53\}.&
\end{array}$$}
\noindent49.
$(n,h,r)=(113,11,3)$:
{\scriptsize$$\begin{array}{llllllll}
\{0, 21, 46\}, &\{0, 18, 41\}, &\{0, 38, 42\}, &\{0, 30, 47\}, &\{0, 20, 51\}, &\{0, 35, 49\}, &\{0, 29, 48, 56\}, \\
\{0, 2, 36\},  &\{0, 11, 37\}, &\{0, 16, 40\}, &\{0, 13, 45, 52\}, &\{0, 1, 44, 54\}, &\{0, 5, 33, 55\}.&&
\end{array}$$}
\noindent50.
$(n,h,r)=(115,7,3)$:
{\scriptsize$$\begin{array}{llllllll}
\hspace{-0.55cm}\{0, 7, 53\}, &\{0, 20, 43\}, &\{0, 34, 42\}, &\{0, 39, 55\},  &\{0, 25, 51, 52\}, &\{0, 36, 41, 54\}, &\{0, 2, 37, 49\},\\
\hspace{-0.55cm}\{0, 10, 40\}, &\{0, 22, 50\}, &\{0, 38, 57\}, &\{0, 14, 31\},  &\{0, 24, 45, 56\}, &\{0, 15, 44, 48\}.&
\end{array}$$}
\noindent51.
$(n,h,r)=(115,19,3)$:
{\scriptsize$$\begin{array}{llllllll}
\{0, 26, 57\}, &\{0, 22, 56\}, &\{0, 14, 54\}, &\{0, 2, 50\}, &\{0, 11, 55\}, &\{0, 5, 51\}, &\{0, 8, 47\}, &\{0, 37, 41\}, \\
\{0, 28, 38\}, &\{0, 25, 32\}, &\{0, 16, 45\}, &\{0, 30, 49\}, &\{0, 42, 43\}, &\{0, 17, 52\}, &\{0, 20, 53\}, &\{0, 13, 36\}.
\end{array}$$}
\noindent52.
$(n,h,r)=(117,9,3)$:
{\scriptsize$$\begin{array}{llllllll}
\hspace{-0.8cm}\{0, 2, 45\}, &\{0, 5, 44\}, &\{0, 13, 47\}, &\{0, 19, 55\}, &\{0, 24, 41\}, &\{0, 32, 54\}, &\{0, 4, 56, 57\}, \\
\hspace{-0.8cm}\{0, 20, 50\}, &\{0, 29, 37\}, &\{0, 35, 51\}, &\{0, 23, 49\}, &\{0, 28, 42\}, &\{0, 33, 40, 58\}, &\{0, 10, 21, 48\}, \\
\hspace{-0.8cm}\{0, 46, 31\}.
\end{array}$$}
\noindent53.
$(n,h,r)=(117,21,3)$:
{\scriptsize$$\begin{array}{llllllll}
\{0, 42, 58\}, &\{0, 11, 56\}, &\{0, 19, 55\}, &\{0, 2, 52\}, &\{0, 5, 53\}, &\{0, 23, 49\}, &\{0, 4, 35\}, &\{0, 17, 37\},  \\
\{0, 38, 46\}, &\{0, 28, 57\}, &\{0, 43, 44\}, &\{0, 33, 40\}, &\{0, 32, 54\}, &\{0, 25, 39\}, &\{0, 41, 51\}, &\{0, 34, 47\}.
\end{array}$$}
\noindent54.
$(n,h,r)=(119,5,3)$:
{\scriptsize$$\begin{array}{llllllll}
\hspace{-0.6cm}\{0, 25, 54\}, &\{0, 39, 57\}, &\{0, 37, 42\}, &\{0, 2, 46\}, &\{0, 43, 53\}, &\{0, 30, 38, 58\}, &\{0, 40, 41, 56\},\\
\hspace{-0.6cm}\{0, 33, 55\}, &\{0, 24, 36\}, &\{0, 45, 52\}, &\{0, 31, 50\},  &\{0, 4, 17, 51\}, &\{0, 14, 23, 49\}, &\{0, 11, 32, 59\}.
\end{array}$$}
\noindent55.
$(n,h,r)=(121,25,3)$:
{\scriptsize$$\begin{array}{llllllll}
\hspace{-0.05cm}\{0, 5, 59\}, &\{0, 2, 58\}, &\{0, 48, 55\}, &\{0, 43, 44\}, &\{0, 42, 52\}, &\{0, 13, 50\}, &\{0, 19, 60\}, &\{0, 4, 38\},\\
\hspace{-0.05cm}\{0, 22, 47\}, &\{0, 45, 53\}, &\{0, 11, 39\}, &\{0, 35, 51\}, &\{0, 31, 57\}, &\{0, 17, 40\}, &\{0, 29, 49\}, &\{0, 32, 46\}.
\end{array}$$}
\noindent56.
$(n,h,r)=(123,9,3)$:
{\scriptsize$$\begin{array}{llllllll}
\hspace{-0.1cm}\{0, 7, 57\}, &\{0, 20, 52\}, &\{0, 5, 49\}, &\{0, 25, 48\},  &\{0, 11, 39, 56\}, &\{0, 41, 51, 59\}, &\{0, 4, 37, 58\}, \\
\hspace{-0.1cm}\{0, 26, 61\}, &\{0, 19, 53\}, &\{0, 14, 43\}, &\{0, 15, 42, 55\}, &\{0, 1, 31, 47\},  &\{0, 22, 24, 60\}.
\end{array}$$}
\noindent57.
$(n,h,r)=(123,21,3)$:
{\scriptsize$$\begin{array}{llllllll}
\{0, 39, 61\}, &\{0, 55, 59\}, &\{0, 20, 56\}, &\{0, 53, 58\}, &\{0, 33, 52\}, &\{0, 42, 50\}, &\{0, 1, 44\}, &\{0, 11, 57\},\\
\{0, 17, 49\}, &\{0, 40, 47\}, &\{0, 14, 48\}, &\{0, 35, 45\}, &\{0, 37, 60\}, &\{0, 29, 31\}, &\{0, 25, 51\}, &\{0, 38, 54\}, \\
\{0, 28, 41\}.&&&&&&&
\end{array}$$}
\noindent58.
$(n,h,r)=(127,25,3)$:
{\scriptsize$$\begin{array}{llllllll}
\hspace{-0.3cm}\{0, 22, 62\}, &\{0, 26, 63\}, &\{0, 14, 61\}, &\{0, 7, 60\}, &\{0, 41, 58\}, &\{0, 46, 56\}, &\{0, 5, 50\}, &\{0, 38, 49\}, \\
\hspace{-0.3cm}\{0, 28, 57\}, &\{0, 34, 59\}, &\{0, 42, 44\}, &\{0, 1, 55\}, &\{0, 4, 52\}, &\{0, 8, 39\}, &\{0, 35, 51\}, &\{0, 19, 32\}, \\
\hspace{-0.3cm}\{0, 23, 43\}.&&&&&&&
\end{array}$$}
\noindent59.
$(n,h,r)=(137,5,3)$:
{\scriptsize$$\begin{array}{llllllll}
\hspace{-0.25cm}\{0, 11, 37\}, &\{0, 42, 50\}, &\{0, 39, 54\}, &\{0, 25, 43\},  &\{0, 5, 51, 68\}, &\{0, 30, 52, 53\}, &\{0, 31, 35, 67\}, \\
\hspace{-0.25cm}\{0, 24, 64\}, &\{0, 33, 60\},&\{0, 20, 49\},  &\{0, 21, 34, 62\}, &\{0, 2, 16, 61\}, &\{0, 9, 47, 66\}, &\{ 0, 7, 55, 65\}, \\
\hspace{-0.25cm}\{0, 12, 56\}.&&&&
\end{array}$$}
\noindent60.
$(n,h,r)=(141,9,3)$:
{\scriptsize$$\begin{array}{llllllll}
\hspace{-0.9cm}\{0, 7, 56\}, &\{0, 19, 52\}, &\{0, 41, 62\}, &\{0, 44, 59\}, &\{0, 31, 55\},  &\{0, 2, 29, 68\},&\{0, 14, 54, 65\}, \\
\hspace{-0.9cm}\{0, 34, 57\}, &\{0, 25, 61\}, &\{0, 32, 60\}, &\{0, 31, 54\}, &\{0, 46, 50\}, &\{0, 35, 48, 53\},&\{0, 47, 63, 64\},\\
\hspace{-0.9cm}\{0, 34, 59\}, &\{0, 52, 62\}, &\{0, 19, 45\}, &\{0, 43, 69\},  &\{0, 38, 58\},  &\{0, 22, 30, 67\}, &\{0, 28, 60, 70\}.
\end{array}$$}
\noindent61.
$(n,h,r)=(141,21,3)$:
{\scriptsize$$\begin{array}{llllllll}
\hspace{0.75cm}\{0, 35, 64\}, &\{0, 14, 56\}, &\{0, 38, 57\}, &\{0, 41, 67\}, &\{0, 43, 65\}, &\{0, 32, 68\},  &\{0, 1, 62, 70\},&\{0, 5, 59, 63\},\\
\hspace{0.75cm}\{0, 7, 47\}, &\{0, 46, 48\}, &\{0, 28, 53\}, &\{0, 31, 51\}, &\{0, 13, 52\}, &\{0, 37, 60\}, &\{0, 11, 45, 55\}, &\{0, 17, 50, 66\}.
\end{array}$$}
\noindent62.
$(n,h,r)=(143,5,3)$:
{\scriptsize$$\begin{array}{llllllll}
\hspace{-0.2cm}\{0, 5, 46\},  &\{0, 9, 44\}, &\{0, 33, 55\}, &\{0, 28, 64\},&\{0, 45, 65, 66\},&\{0, 11, 43, 61\}, &\{0, 34, 58, 60\},\\
\hspace{-0.2cm}\{0, 10, 47\}, &\{0, 51, 63\}, &\{0, 31, 38\}, &\{0, 4, 23, 71\}, &\{0, 14, 53, 68\}, &\{0, 8, 57, 70\}, &\{0, 27, 52, 69\},\\
\hspace{-0.2cm}\{0, 40, 56\}, &\{0, 29, 59\}.&&
\end{array}$$}
\noindent63.
$(n,h,r)=(145,25,3)$:
{\scriptsize$$\begin{array}{llllllll}
\{0, 8, 72\}, &\{0, 39, 71\}, &\{0, 25, 70\}, &\{0, 19, 68\}, &\{0, 4, 69\}, &\{0, 56, 67\}, &\{0, 5, 66\}, &\{0, 28, 59\},\\
\{0, 7, 50\}, &\{0, 46, 62\}, &\{0, 26, 48\}, &\{0, 34, 63\}, &\{0, 47, 57\}, &\{0, 42, 44\}, &\{0, 23, 60\}, &\{0, 17, 58\},\\
\{0, 20, 55\}, &\{0, 14, 54\}, &\{0, 1, 53\}, &\{0, 38, 51\}.&&&&
\end{array}$$}
\noindent64.
$(n,h,r)=(147,9,3)$:
{\scriptsize$$\begin{array}{llllllll}
\hspace{-0.55cm}\{0, 7, 58\}, &\{0, 11, 47\}, &\{0, 25, 55\}, &\{0, 39, 53\}, &\{0, 40, 57\}, &\{0, 10, 66, 70\}, &\{0, 24, 26, 61\},\\
\hspace{-0.55cm}\{0, 28, 71\}, &\{0, 29, 73\}, &\{0, 49, 50\}, &\{0, 27, 48\}, &\{0, 5, 64, 72\}, &\{0, 15, 31, 69\},  &\{0, 19, 32, 65\},\\
\hspace{-0.55cm}\{0, 18, 52\},  &\{0, 20, 62\},&\{0, 23, 68\},  &\{0, 41, 63\}.&&&
\end{array}$$}
\noindent65.
$(n,h,r)=(151,25,3)$:
{\scriptsize$$\begin{array}{llllllll}
\hspace{-0.15cm}\{0, 35, 75\}, &\{0, 17, 74\}, &\{0, 19, 73\}, &\{0, 48, 68\}, &\{0, 56, 70\}, &\{0, 71, 72\}, &\{0, 4, 64\}, &\{0, 13, 59\}, \\
\hspace{-0.15cm}\{0, 7, 41\}, &\{0, 23, 62\}, &\{0, 28, 53\}, &\{0, 29, 67\}, &\{0, 31, 63\}, &\{0, 42, 58\}, &\{0, 49, 51\}, &\{0, 10, 65\}, \\
\hspace{-0.15cm}\{0, 43, 69\}, &\{0, 44, 66\}, &\{0, 50, 61\}, &\{0, 47, 52\}, &\{0, 37, 45\}.&&&
\end{array}$$}
\noindent66.
$(n,h,r)=(165,9,3)$:
{\scriptsize$$\begin{array}{llllllll}
\hspace{0.35cm}\{0, 13, 80\}, &\{0, 61, 66\}, &\{0, 24, 63\}, &\{0, 19, 52, 69\}, &\{0, 25, 48, 76\}, &\{0, 14, 29, 70\}, &\{0, 30, 34, 77\}, \\
\hspace{0.35cm}\{0, 49, 31\}, &\{0, 45, 71\}, &\{0, 53, 55\}, &\{0, 16, 60, 81\}, &\{0, 7, 42, 82\}, &\{0, 1, 38, 74\}, &\{ 0, 10, 64, 72\}, \\
\hspace{0.35cm}\{0, 58, 78\}, &\{0, 32, 59\},  &\{0, 22, 68, 79\}.
\end{array}$$}
\noindent67.
$(n,h,r)=(169,25,3)$:
{\scriptsize$$\begin{array}{llllllll}
\{0, 7, 62\}, &\{0, 10, 84\}, &\{0, 19, 69\}, &\{0, 77, 82\}, &\{0, 46, 78\}, &\{0, 80, 81\}, &\{0, 14, 70\}, &\{0, 58, 75\}, \\
\{0, 38, 66\}, &\{0, 4, 67\}, &\{0, 11, 64\}, &\{0, 23, 68\}, &\{0, 42, 73\}, &\{0, 34, 71\}, &\{0, 41, 76\}, &\{0, 59,79\}, \\
\{0, 44, 60\}, &\{0, 49, 57\},  &\{0, 47, 72\}, &\{0, 22, 51\}, &\{0, 43, 83\}, &\{0, 48, 61\}, &\{0, 39, 65\}, &\{0, 52, 54\}.
\end{array}$$}
\noindent68.
$(n,h,r)=(171,9,3)$:
{\scriptsize$$\begin{array}{llllllll}
\hspace{0.05cm}\{0, 1, 77\}, &\{0, 35, 83\}, &\{0, 2, 66\}, &\{0, 4, 62\},  &\{0, 10, 17, 71\}, &\{0, 28, 49, 79\}, &\{ 0, 20, 44, 80\},\\
\hspace{0.05cm}\{0, 16, 68\}, &\{0, 33, 75\}, &\{0, 46, 85\},  &\{0, 55, 69\},  &\{0, 50, 73, 84\}, &\{0, 26, 57, 82\}, &\{0, 19, 32, 72\},\\
\hspace{0.05cm}\{0, 65, 70\},  &\{0, 8, 67\}, &\{0, 38, 81\}, &\{0, 27, 45, 74\}, &\{0, 15, 37, 78\}.
\end{array}$$}
\noindent69.
$(n,h,r)=(189,9,3)$:
{\scriptsize$$\begin{array}{llllllll}
\hspace{-0.45cm}\{0, 1, 76\}, &\{0, 2, 82\}, &\{0, 43, 69\}, &\{0, 48, 78\}, &\{0, 28, 74, 92\}, &\{0, 54, 71, 94\}, &\{0, 8, 70, 85\},  \\
\hspace{-0.45cm}\{0, 33, 67\}, &\{0, 16, 57\}, &\{0, 11, 61\}, &\{0, 5, 56\}, &\{0, 10, 35, 93\}, &\{0, 49, 81, 88\}, &\{0, 19, 55, 84\},\\
\hspace{-0.45cm}\{0, 44, 89\},  &\{0, 73, 87\}, &\{0, 21, 68\}, &\{0, 20, 72\}, &\{0, 38, 60, 91\},  &\{0, 27, 86, 90\}, &\{0, 42, 66, 79\}.
\end{array}$$}
\noindent70.
$(n,h,r)=(195,9,3)$:
{\scriptsize$$\begin{array}{llllllll}
\hspace{-0.25cm}\{0, 45, 93\}, &\{0, 29, 72\}, &\{0, 25, 88\}, &\{0, 28, 84\},&\{0, 36, 55, 90\}, &\{0, 7, 15, 77\}, &\{0, 21, 71, 89\},\\
\hspace{-0.25cm}\{0, 46, 85\}, &\{0, 53, 69\}, &\{0, 23, 80\},   &\{0, 40, 81\}, &\{0, 22, 60, 73\}, &\{0, 30, 67, 94\},  &\{0, 31, 32, 97\},\\
\hspace{-0.25cm}\{0, 44, 86\}, &\{0, 83, 87\}, &\{0, 47, 58\}, &\{0, 2, 26, 78\}, &\{0,5,79,96\}, &\{0, 61, 75, 95\}, &\{0, 10, 59, 92\}.
\end{array}$$}
\noindent71.
$(n,h,r)=(213,9,3)$:
{\scriptsize$$\begin{array}{llllllll}
\hspace{0.75cm}\{0, 5,104\}, &\{0, 15, 94\}, &\{0, 34, 95\},  &\{0, 19, 71\},  &\{0, 62, 64, 106\}, &\{0, 56, 84, 92\}, &\{0, 4, 86, 97\}, \\
\hspace{0.75cm}\{0, 45, 78\}, &\{0, 49, 69\}, &\{ 0, 55, 103\},  &\{0, 60, 70, 87\}, &\{0, 26, 57, 98\}, &\{0, 46, 59, 96\}, &\{0, 1, 68, 91\}, \\
\hspace{0.75cm}\{0, 74, 81\},  &\{0, 53, 85\}, &\{0, 16, 89\}, &\{0, 25, 43, 83\}, &\{0, 38, 77, 101\}, &\{0, 14, 80, 102\},  &\{0, 30, 51,105\},\\
\hspace{0.75cm}\{0, 29, 76\},  &\{0, 65, 100\}.
\end{array}$$}
\noindent72.
$(n,h,r)=(219,9,3)$:
{\scriptsize$$\begin{array}{llllllll}
\hspace{0.55cm}\{0, 2, 88\}, &\{0, 42,105\},&\{0, 1, 82\}, &\{0, 15, 99\}, &\{0, 11, 62\},  &\{0, 14, 21, 108\}, &\{0, 24, 54,98\},\\
\hspace{0.55cm}\{0, 28, 92\}, &\{0, 41, 91\}, &\{0, 5, 76\},  &\{0, 23,103\},  &\{0, 26, 61, 93\},  &\{0, 43, 83, 100\}, &\{0, 18, 97, 107\},\\
\hspace{0.55cm}\{0, 53, 78\}, &\{0, 19, 52\}, &\{0, 48, 95\}, &\{0, 22, 49\},  &\{0, 4, 59, 72\}, &\{0, 46, 77, 106\}, &\{0, 34, 73, 109\},\\
\hspace{0.55cm}\{0, 20, 90\}, &\{0, 16, 85\}, &\{0, 37, 102\},  &\{0, 56, 101\},   &\{0, 38,96, 104\}.
\end{array}$$}
\noindent73.
$(n,h,r)=(237,9,3)$:
{\scriptsize$$\begin{array}{llllllll}
\hspace{-0.9cm}\{0, 1, 91\}, &\{0,25,88\}, &\{0, 2, 82\}, &\{0, 14, 108\}, &\{0,  22, 37, 101\}, &\{ 0, 32, 51, 100\},  \\
\hspace{-0.9cm}\{0, 36, 89\},  &\{0, 47, 81\}, &\{0, 40, 114\},  &\{0, 10, 60, 117\},  &\{0, 55,111,116\}, &\{0, 71, 84, 104\},    \\
\hspace{-0.9cm}\{0, 58, 99\},  &\{0, 83, 106\}, &\{0, 35, 97\}, &\{0, 38, 67, 115\}, &\{0, 17, 95, 103\}, &\{0, 26, 102, 118\}, \\
\hspace{-0.9cm}\{0, 52, 70\}, &\{0, 24, 93\},  &\{0, 44, 87\},  &\{0, 4, 11, 109\}, &\{0, 46,73, 112\}, &\{0, 54, 85, 113\},\\
\hspace{-0.9cm}\{0, 21, 96\},  &\{0, 65, 110\},  &\{0, 30, 72\}.
\end{array}$$}
\noindent74.
$(n,h,r)=(243,9,3)$:
{\scriptsize$$\begin{array}{llllllll}
\hspace{-0.9cm}\{0, 4, 84\}, &\{0, 70, 119\},  &\{0, 8, 105\}, &\{0, 16, 111, 118\}, &\{0, 35, 50, 113\},  &\{ 0, 13, 40, 116\},  \\
\hspace{-0.9cm}\{0, 5, 71\},&\{0, 62, 115\}, &\{0, 51, 107\}, &\{0, 17, 18, 117\}, &\{0, 25, 39, 98\}, &\{0, 29, 93, 121\},\\
\hspace{-0.9cm}\{0, 33, 85\},  &\{0, 57, 87\}, &\{0, 47,79\},  &\{0, 48, 109, 120\}, &\{0,  20, 89, 110\}, &\{0, 2, 77, 114\}, \\
\hspace{-0.9cm}\{0, 42, 86\}, &\{0, 31, 91\}, &\{0, 26, 81\},&\{0, 58, 94,104\}, &\{0, 41, 65, 108\},  &\{0, 23, 68, 106\},\\
\hspace{-0.9cm}\{0, 82, 101\}, &\{0, 74,96\},  &\{0, 54,88\}.
\end{array}$$}
\noindent75.
$(n,h,r)=(261,9,3)$:
{\scriptsize$$\begin{array}{llllllll}
\hspace{-0.25cm}\{0, 8, 90\}, &\{0, 7, 80 \}, &\{0, 1, 101\}, &\{0, 21, 26, 123\},&\{0, 54, 89, 128\}, &\{0, 14, 69, 130\},\\
\hspace{-0.25cm}\{0, 87, 110\}, &\{0, 58, 129\}, &\{0, 70, 111\},  &\{0, 63, 78, 122\}, &\{0,47, 99, 124\},  &\{0, 24, 88, 119\},\\
\hspace{-0.25cm}\{0, 65, 85\},  &\{0, 43, 103\},  &\{0, 27,120\}, &\{0, 48, 86, 104\}, &\{0, 40, 57, 106\},  &\{0, 32, 62, 115\},\\
\hspace{-0.25cm}\{0, 34, 109\},  &\{0, 33, 125\}, &\{0, 79, 108\}, &\{0, 42, 46, 114\}, &\{0, 105, 107, 118\}, &\{0, 36, 117, 127\},\\
\hspace{-0.25cm}\{0, 37, 121\}, &\{0, 19, 113\}, &\{0, 28, 50, 126\}, &\{0, 45,96, 112\}.
\end{array}$$}
\noindent76.
$(n,h,r)=(267,9,3)$:
{\scriptsize$$\begin{array}{llllllll}
\hspace{-0.8cm}\{0, 2, 105\}, &\{0, 7, 109\}, &\{0, 4, 88\}, &\{0, 8, 106\},  &\{0, 74, 91, 123\}, &\{0, 55, 96, 125\}, \\
\hspace{-0.8cm}\{0, 33, 100\}, &\{0, 25, 97\}, &\{0, 51, 79\}, &\{0, 44, 90\}, &\{0, 24, 93, 104\}, &\{0, 34, 57, 117\}, \\
\hspace{-0.8cm}\{0, 87, 122\}, &\{0, 52, 118\}, &\{0, 27, 75\},   &\{0, 10, 40, 129\}, &\{0, 19, 131, 132\}, &\{0, 5, 76, 121\}, \\
\hspace{-0.8cm}\{0, 18, 95\}, &\{0, 111, 126\}, &\{0, 13, 86\}, &\{0, 43, 82, 124\}, &\{0, 63, 101, 127\}, &\{0, 78, 92, 128\},  \\
\hspace{-0.8cm}\{0, 65, 133\}, &\{0, 61, 108\}, &\{0, 56, 115\},   &\{0, 62, 99, 120\},  &\{0, 16, 110, 130\}, &\{0, 54, 85,107\}.
\end{array}$$}
\noindent77.
$(n,h,r)=(285,9,3)$:
{\scriptsize$$\begin{array}{llllllll}
\hspace{-1.2cm}\{0,16, 113\}, &\{0, 30,123\}, &\{0, 24, 107\}, &\{0, 62, 111\},&\{0, 46, 128, 142\}, &\{0, 63, 92, 100\}, \\
\hspace{-1.2cm}\{0, 72, 137\}, &\{0, 74, 127\}, &\{0, 69, 101\},&\{0, 40, 84, 139\},   &\{0, 47, 51, 124\},&\{0, 35, 45, 133\},\\
\hspace{-1.2cm}\{0, 70, 95\}, &\{0, 58, 125\}, &\{0, 91, 130\},  &\{0, 5, 18, 134\}, &\{0, 79, 81, 138\},   &\{0, 71, 78, 114\},\\
\hspace{-1.2cm}\{0, 75, 136\},&\{0, 27, 103\}, &\{0, 86, 119\},  &\{0, 26, 94, 115\}, &\{0, 19, 121,141\}, &\{0, 1, 110, 132\},  \\
\hspace{-1.2cm}\{0, 60, 112\},  &\{0, 64, 120\},&\{0, 85, 126\},  &\{0, 23, 34, 140\},  &\{0, 28, 66, 108\}, &\{0, 17, 48, 135\}, \\
\hspace{-1.2cm}\{0, 90, 105\},  &\{0, 54, 104\}.
\end{array}$$}
\noindent78.
$(n,h,r)=(291,9,3)$:
{\scriptsize$$\begin{array}{llllllll}
\hspace{-0.2cm}\{0, 5, 101\}, &\{0, 18, 100\}, &\{0, 30, 86\}, &\{0, 32, 121\}, &\{0,61,97,145\},&\{0, 10, 44, 129\},\\
\hspace{-0.2cm}\{0, 42, 73\}, &\{0, 66, 131\}, &\{0, 68, 144\}, &\{0, 53, 116\},   &\{0, 21, 70, 141\},  &\{0, 11, 25, 133\},\\
\hspace{-0.2cm}\{0, 54, 104\},  &\{0, 23, 117\},&\{0, 1, 92, 137\}, &\{0, 7, 118, 135\}, &\{0, 29, 112, 132\},  &\{0,  59, 98, 138\},\\
\hspace{-0.2cm}\{0, 26, 107\}, &\{0, 35, 113\}, &\{0, 33, 46, 123\}, &\{0, 47, 88, 140\}, &\{0, 62, 134, 142\}, &\{0, 4, 99, 114\},\\
\hspace{-0.2cm}\{0, 28, 143\}, &\{0, 16, 74, 125\}, &\{0,  22, 60, 127\},  &\{0, 75, 102, 139\}, &\{0, 69, 124, 126\}, &\{0, 24, 43, 130\}.
\end{array}$$}
\noindent79.
$(n,h,r)=(309,9,3)$:
{\scriptsize$$\begin{array}{llllllll}
\hspace{-0.8cm}\{0, 1, 122\}, &\{0, 20, 149\}, &\{0, 39, 130\}, &\{0, 7, 26, 154\}, &\{0, 5, 60, 132\}, &\{0, 28, 36,120\},\\
\hspace{-0.8cm}\{0, 64, 134\},  &\{0, 89, 123\}, &\{0, 37, 113\}, &\{0, 52, 81, 146\},  &\{ 0, 48, 99,153\},&\{0, 23,53,148\},  \\
\hspace{-0.8cm}\{0, 96, 143\}, &\{0, 67, 111\}, &\{0, 46, 126\},    &\{0, 27,49,136\}, &\{0, 62, 83, 140\}, &\{0, 14, 104, 145\},\\
\hspace{-0.8cm}\{0, 33, 107\}, &\{0, 85, 135\}, &\{0, 17, 110\},   &\{0, 35, 138, 151\}, &\{0, 45, 61, 124\},&\{0, 86, 118,142\},  \\
\hspace{-0.8cm}\{0, 10, 98\}, &\{0, 77, 150\}, &\{0, 42, 108\},  &\{0, 18, 133, 137\},&\{0, 71, 114, 139\},    &\{ 0, 11, 112, 152\},\\
\hspace{-0.8cm}\{0, 15, 97\}, &\{0, 59, 117\},  &\{0, 100, 102\},&\{0, 75, 106, 144\}.
\end{array}$$}

\section{$(n,\{3,4,5\},1)$-PDPs with leave $[h]^r$}\label{app-D}
The base blocks of an $(n,\{3,4,5\},1)$-PDP with leave $[h]^r$ are listed below:\\

\noindent1.
$(n,h,r)=(31,5,4)$:
{\scriptsize\{0, 6, 13\}, \{0, 3, 5, 14, 15\}.}\\

\noindent2.
$(n,h,r)=(37,5,5)$:
{\scriptsize\{0, 9, 17\}, \{0, 6, 13\}, \{0, 2, 3, 14, 18\}.}\\

\noindent3.
$(n,h,r)=(43,5,5)$:
{\scriptsize\{0, 18, 11\}, \{0, 12, 21\}, \{0, 6, 14\}, \{0, 4, 17, 19, 20\}.}\\

\noindent4.
$(n,h,r)=(49,5,4)$:
{\scriptsize\{0, 13, 24\}, \{0, 10, 22\}, \{0, 19, 21\}, \{0, 7, 16\}, \{0, 3, 17, 18, 23\}.}\\

\noindent5.
$(n,h,r)=(55,5,4)$:
{\scriptsize\{0, 5, 27\}, \{0, 19, 20\}, \{0, 18, 25\}, \{0, 10, 26\}, \{0, 13, 24\}, \{0, 6, 9, 21, 23\}.}\\

\end{CJK*}


\begin{thebibliography}{99}

\bibitem{Difference-families}
Abel R J R and Buratti M.  ``Difference families", in Handbook of Combinatorial Designs, 2nd edn, Colbourn C J and Dinitz J H, (Editors), Chapman and Hall, CRC, Boca Raton, FL, 2007: 392-410.

\bibitem{MGDD}
Abel R J R and Assaf A M. Modified group divisible designs with block size $5$ and $\lambda=1$. Discrete Math., 2002, 256(1-2): 1-22.

\bibitem{5-MGDD}
Abel R J R and Assaf A M. Modified group divisible designs with block size $5$ and even index. Discrete Math., 2008, 308(15): 3335-3351.


\bibitem{Abrham}
Abrham J. Perfect systems of difference sets-A survey. Ars Combin., 1984, 17A: 5-36.

\bibitem{3-MGDD}
Assaf A. Modified group divisible designs. Ars Combin., 1990, 29: 13-20.


\bibitem{Bermond}
Bermond J C, Kotzig A and Turgeon J. On a combinatorial problem of antennas in radio astronomy. In: Proc. 18th Hungarian Combinatorial Colloquium, North Holland, 1976: 135-149.

\bibitem{Beth}
Beth T, Jungnickel D and Lenz H. Design theory. Vol. I. Cambridge University Press Cambridge, 1999.

\bibitem{caoh}
Cao H, Wang L and Wei R. The existence of HGDDs with block size four and its application to double frames. Discrete Math., 2009, 309: 945-949.

\bibitem{Chang}
Chang Y and Miao Y. Constructions for optimal optical orthogonal codes. Discrete Math., 2003, 261(1-3): 127-139.

\bibitem{cklh}
Chee Y M, Kiah H M, Ling S and Wei H. Geometric orthogonal codes of size larger than optical orthogonal codes. IEEE Trans. Inf. Theory, 2017, 64(4): 2883-2895.

\bibitem{ChenZ}
Chen Z. The existence of balanced difference families and perfect difference families. Master Degree Thesis, Guangxi Normal University, 2008.

\bibitem{Doty}
Doty D and Winslow A. Design of geometric molecular bonds. IEEE Trans. Mol. Biol. Multi-Scale Commun., 2017, 3(1): 13-23.

\bibitem{Feng}
Feng T, Wang X and Chang Y. Semi-cyclic holey group divisible designs with block size three. Des. Codes Cryptogr., 2015, 74(2): 301-324.

\bibitem{Radar}
Ge G, Ling A and Ying M. A systematic construction for radar arrays. IEEE Trans. Inf. Theory, 2008, 54(1): 410-414.

\bibitem{GeG}
Ge G, Miao Y and Sun X. Perfect difference families, perfect difference matrices, and related combinatorial structures. J. Combin. Des., 2010, 18(6): 415-449.

\bibitem{Golomb}
Golomb S W and Taylor H. Two-dimensional synchronization patterns for minimum ambiguity. IEEE Trans. Inf. Theory, 1982, 28(4): 600-604.

\bibitem{Huang}
Huang J H and Skiena S S. Gracefully labeling prisms. Ars Combin., 1994, 38: 225-236.

\bibitem{Jiang}
Jiang J, Wu D and Fan P. General constructions of optimal variable-weight optical orthogonal codes. IEEE Trans. Inf. Theory, 2011, 7: 4488-4496.

\bibitem{Laufer}
Laufer P J and Turgeon J M. On a conjecture of Paul Erd\"os for perfect systems of difference sets. Discrete Math., 1983, 47: 255-266.

\bibitem{Rothemund}
Rothemund P. Using lateral capillary forces to compute by self-assembly. Proc. National Academy of Sciences, 2000, 97(3): 984-989.

\bibitem{Shearer}
Shearer J B. ``Difference triangle sets", in Handbook of Combinatorial Designs, 2nd edn,
Colbourn C J and Dinitz J H (Editors), Chapman and Hall, CRC, Boca Raton, FL, 2007: 436-440.

\bibitem{LD}
Simpson J E. Langford sequences: Perfect and hooked. Discrete Math., 1983, 44(1): 97-104.

\bibitem{WLD}
Wang L, Cai L, Feng T, Tian Z and Wang X. Geometric orthogonal codes and generalized perfect difference families. Submitted.

\bibitem{WuD}
Wu D, Cheng M and Chen Z. Perfect difference families and related variable-weight optical orthogonal codes. Australas. J. Combin., 2013, (55): 153-166.

\bibitem{yang}
Yang G C. Variable-weight optical orthogonal codes for CDMA networks with multiple performance requirements. IEEE Trans. Commun., 1996, 44: 47-55.

\bibitem{Zhang}
Zhang Z and Tu C. New bounds for the sizes of radar arrays. IEEE Trans. Inf. Theory, 1994, 40: 1672-1678.



\end{thebibliography}
\end{document}